\RequirePackage{fix-cm}
\documentclass[smallextended]{svjour3}       
\smartqed  
\usepackage{graphicx}
\usepackage{amsmath,amssymb}
\usepackage{subfigure}
\usepackage[percent]{overpic}

\usepackage{color}

\newcommand{\R}{\mathbb{R}}
\newcommand{\ds}{\displaystyle}

\spnewtheorem{assumption}{Assumption}{\bf}{\it}
\spnewtheorem{rem}{Remark}{\bf}{\it}
\spnewtheorem{alg}{Algorithm}{\bf}{\rm}

\begin{document}

\title{Coupling  regularization and adaptive $hp$-BEM for the solution of  
a delamination problem
}

\titlerunning{$hp$-adaptive BEM for regularized hemivariational inequalities}        

\author{Nina~Ovcharova         \and
        Lothar~Banz 
}


\institute{N. Ovcharova \at
               Universit\"at der Bundeswehr M\"unchen, D-85577 Neubiberg/Munich, Germany \\
              \email{nina.ovcharova@unibw.de}           
           \and
           L. Banz \at
              Department of Mathematics, University of Salzburg, Hellbrunner Stra{\ss}e 34, 5020 Salzburg, Austria\\
            \email{lothar.banz@sbg.ac.at} 
}

\date{Received: date / Accepted: date}

\maketitle

\begin{abstract}
In this paper, we couple  regularization techniques 
with the adaptive $hp$-version of the boundary element method ($hp$-BEM) for the efficient numerical solution of linear elastic problems with nonmonotone contact boundary conditions. As a model example we treat the delamination of composite structures with a contaminated interface layer.  This problem has a weak formulation in terms of a nonsmooth variational inequality. The resulting  hemivariational inequality (HVI) 
is first regularized and then, discretized by an adaptive $hp$-BEM. We give conditions for the uniqueness of the solution and provide an a-priori error estimate. Furthermore,  we derive an a-posteriori error estimate
for the nonsmooth variational problem based on  a novel regularized mixed formulation, thus enabling $hp$-adaptivity.
Various numerical experiments illustrate the behavior, strengths and weaknesses of the proposed high-order approximation scheme. 
\end{abstract}

\section{Introduction}

The interest in robust and efficient numerical methods for the solution of nonsmooth problems in nonmonoton contact like  adhesion and  delamination problems increases constantly. The nonsmoothness comes from the nonsmooth data of the problems itself, in particular from nonmonotone multivalued physical laws involved in the boundary conditions that lead to nonsmooth functionals in the variational formulation. 
There are several approaches to treat this non-differentiability. We can first 
discretize by finite elements, and then solve the nonconvex optimization problems by novel  nonsmooth optimization methods \cite{Noll_Ovcharova}, or first regularize the nonsmooth functional, then discretize by finite element methods and finally use standard optimization solvers \cite{Ovcharova2012}. Since  the interesting
nonmonotone behavior takes place only in a boundary part of the involved linear elastic  bodies, the domain hemivariational inequaliy can also be formulated as 
a boundary hemivariational inequality 
via the Poincar\'{e}-Steklov operator. 
Therefore, the boundary element methods (BEM) are the methods of choice for discretization. There are number of works exploiting the $h$- and the high-order $hp$-version of the BEM for monotone contact problems. 
Advanced $hp$-BEM discretizations, a-priori and a-posteriori error estimates for unilateral Signorini problems or contact problems with monotone friction are established in \cite{banz2013posteriori,Banz2015,chernov2008hp,MS-1,maischak2007adaptive,stephan2009hp}. However, successfully applied  $hp$-BEM techniques for nonsmooth nonconvex  variational inequalities  are still missing. Pure $h$-BEM  with lowest polynomial degree has been investigated for the first time by Nesemann and Stephan  \cite{Nesemann} for unilateral contact  with adhesion. 
They study the existence and uniqueness, and introduce a residual error indicator.
 
In this paper, we extend the  approximation scheme from \cite{Ovcharova2015} (dealing with  the $h$ - version of the BEM) to the $hp$-version of the BEM for the regularized problem. 
More precisely,  
the nonsmooth variational inequality is first regularized and then discretized by $hp$-adaptive BEM. 
For the Galerkin solution of the regularized problem we provide an a-priori error estimate 
and derive a reliable a-posteriori error estimate based on a {\it novel}  regularized mixed formulation and thus, enabling $hp$-adaptivity. All our theoretical results are illustrated with various numerical experiments for a contact problem with adhesion. 
We also provide conditions for the uniqueness of the solution that sharpen the results due to \cite{Nesemann}.  

\textbf{Notation:} We denote with $C_i$, $c_i$ and such alike generic constants, which can take different values at different positions. We use bold font to indicate vector-valued variables, e.g.~$\mathbf {u} \in \mathbf{H}^1(\Omega):=[H^1(\Omega)]^d$. If there are too many indices, one of them
is written as a superscript, 
e.g.~$\mathbf{u}_{\varepsilon, hp}=\mathbf{u}^{\varepsilon}_{hp}$.

\section{A nonmonotone boundary value problem from delamination}
Let $\Omega \subset \R^d \, (d= 2, 3)$ be a bounded  domain with Lipschitz boundary $\Gamma=\partial \Omega$. We assume that the boundary is decomposed into three disjoint open parts $\Gamma_D, \Gamma_N$,
and $\Gamma_C$ such that $\Gamma = \overline{\Gamma}_D \cup
\overline{\Gamma}_N \cup \overline{\Gamma}_C$ and, moreover, the measures of $\Gamma_C$ and $\Gamma_D$ are positive. We fix an elastic body occupying $\overline{\Omega}$. The body is subject to volume force ${\mathbf f} \in [L^2(\Omega)]^d$ and $g\in H^{1/2}(\Gamma_C)$, $g\geq 0$, is a gap function associating every point $x\in \Gamma_C$ with its distance to the rigid obstacle measured in the direction of the unit outer normal vector $\mathbf{n}(x)$. The body is fixed along $\Gamma_D$, surface tractions 
$\mathbf {t}  \in [L^2(\Gamma_N)]^d$ act on $\Gamma_N$, and  on the part $\Gamma_C$ a nonmonotone, generally multivalued boundary condition holds. 

Further, $\varepsilon (\mathbf {u}) = \frac{1}{2}(\nabla {\mathbf u}  + \nabla {\mathbf u}^T)$ denotes the  linearized strain tensor  and $\sigma ( \mathbf {u}) = \mathcal{C}  \,: \, \varepsilon (\mathbf {u})$ stands for the stress tensor, where $\mathcal{C}$ is the Hooke tensor, assumed to be 
uniformly positive definite with $L^\infty$ coefficients. The boundary stress vector  can be decomposed further into the normal, respectively, the tangential stress: 
\[
\sigma_n= \sigma(\mathbf{u}) {\mathbf n} \cdot {\mathbf n}, \quad \sigma_t=\sigma(\mathbf{u}) {\mathbf n} - \sigma_n {\mathbf n}.
\]

\begin{figure}[b]
\centering
\includegraphics[trim = 0cm 0cm 0cm 0cm, clip,width=0.9\textwidth]{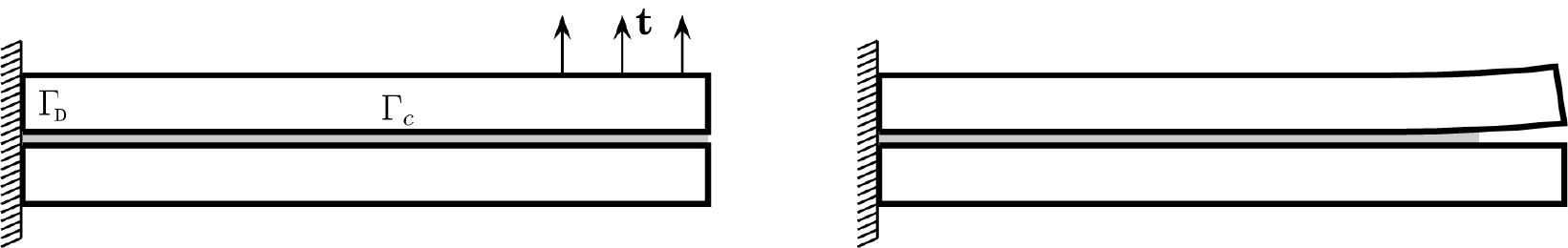}
\caption{Reference configuration for the 2D benchmark under loading. Under applied traction force $\mathbf{t}$ the crack front propagates to the left} \label{delpr}
\end{figure}

By assuming that the structure is symmetric and the forces applied to the upper and lower part are the same, it suffices to consider only the upper half of the specimen represented by $\overline{\Omega}$, see Fig. \ref{delpr} left for the  2D benchmark problem.

\begin{figure}[tbp]
\centering \mbox{
  \subfigure[Load-displacement curve for different contamination concentrations, see \cite{Gudladt}]{
\includegraphics[trim = 0cm 9cm 0cm 0cm, clip,width=0.5\textwidth]{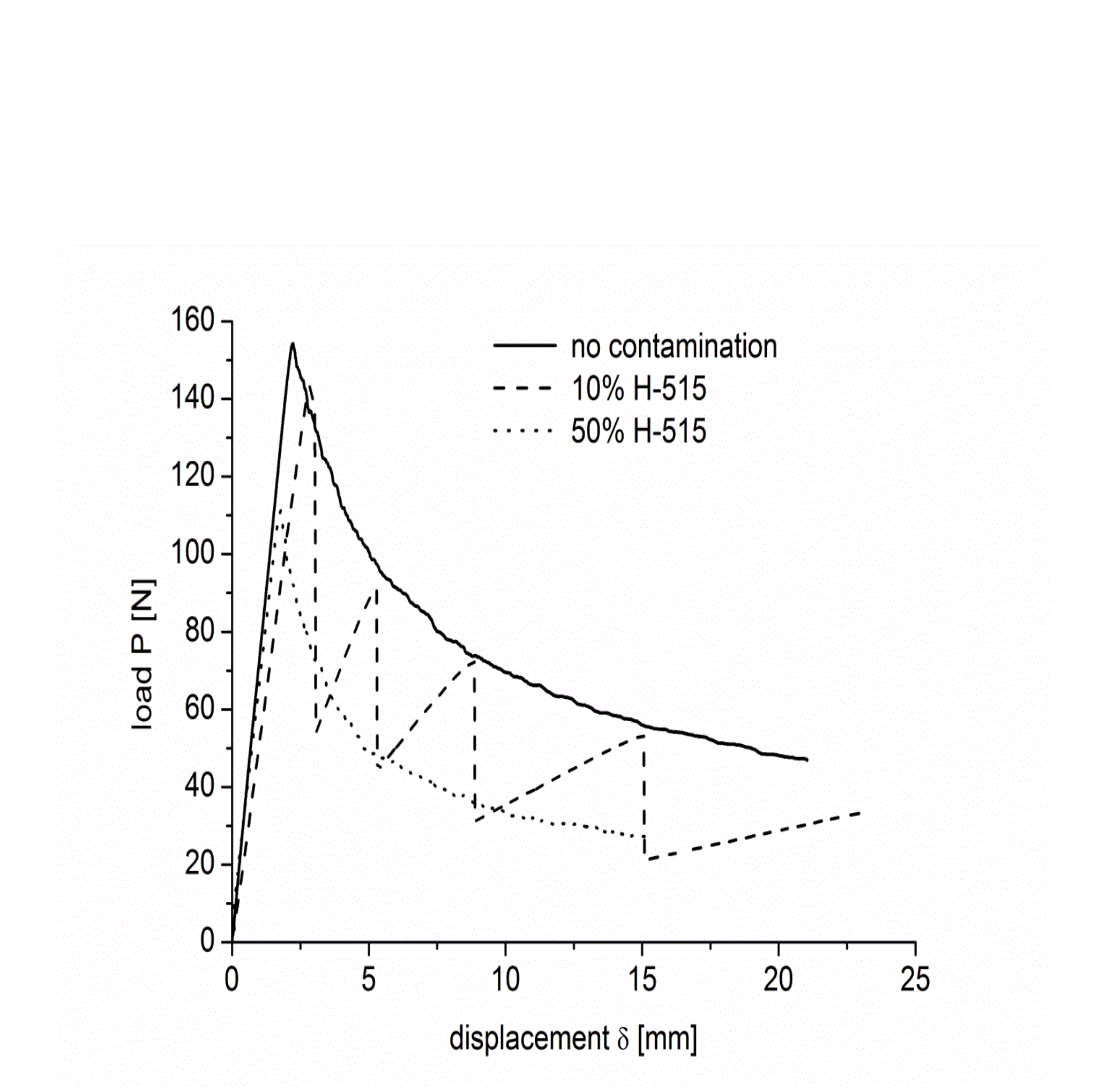}} \ 
\subfigure[A nonmonotone adhesion law]{
	\includegraphics[trim = 50mm 150mm 0mm 0mm, clip,width=65.0mm, keepaspectratio]{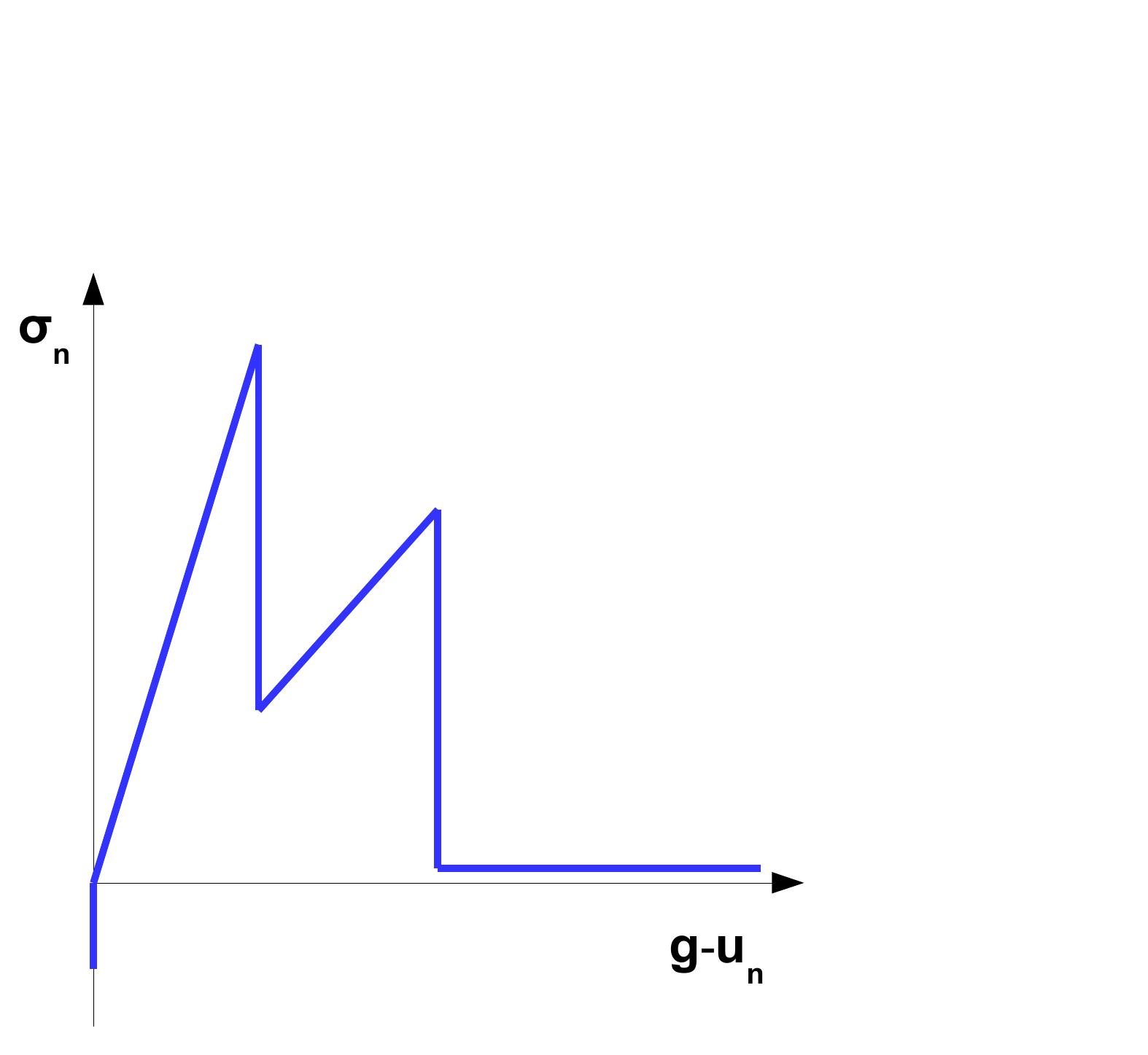}} } 
\caption{Nonmonotone adhesion law} \label{01}
\end{figure}
The delamination problem under consideration is the following: Find  ${\mathbf u} \in \mathbf{H}^1(\Omega):=[H^1(\Omega)]^d$ such that 
\begin{subequations}
\begin{alignat}{2}
-\mbox{div} \; \sigma ({\mathbf u})&= {\mathbf f}  & \quad & \mbox{in} \; \Omega \label{eq1} \\
{\mathbf u}&=0  & \quad & \mbox{on} \; \Gamma_D  \\
\sigma({\mathbf u}){\mathbf n} &= {\mathbf t}  & \quad & \mbox{on}\; \Gamma_N \\
u_n &\leq g  & \quad & \mbox{on}\; \Gamma_C \label{unilateral} \\
\sigma_t({\mathbf u})&=0   & \quad & \mbox{on} \; \Gamma_C\label{eq:notangstress} \\
-\sigma_n ({\mathbf u}) &\in \partial f(u_n) & \quad & \mbox{on} \; \Gamma_C . \label{nbc}
\end{alignat}
\end{subequations}
The contact law (\ref{nbc}), written as a differential inclusion by means  of  the Clarke subdifferential $\partial f$ \cite{Clarke} of a locally Lipschitz function $f$, describes the  nonmonotone, multivalued behavior of the adhesive. More precisely, $\partial f$ is the physical law between the normal component $\sigma_n$ of the stress boundary vector and the normal component $u_n= \mathbf{u} \cdot \mathbf{n}$ of the displacement $\mathbf{u}$ on $\Gamma_C$.  A typical zig-zagged nonmonotone adhesion law is shown in Figure \ref{01}(b).\\
To give a variational formulation of the above boundary value problem we define
\begin{align*}
\mathbf{H}^1_D(\Omega)=\left\{\mathbf{v} \in \mathbf{H}^1(\Omega) \, : \, \mathbf{v}|_{\Gamma_D}=0 \right\}, \quad 
\mathcal{K}=\{ \mathbf{v} \in \mathbf{H}^1_D(\Omega) \, :\,  v_n \leq g \text{ a.e.~on } \Gamma_C\}
\end{align*}
and introduce the $\mathbf{H}^1_D(\Omega)$-coercive and continuous bilinear form of linear elasticity
\[
a(\mathbf{u},\mathbf{v})= \int_{\Omega} \sigma (\mathbf{u}) : \varepsilon (\mathbf{v}) \, dx.
\]
Multiplying the equilibrium equation  (\ref{eq1}) by $\mathbf{v-u}$, integrating over $\Omega$ and applying the divergence theorem yields
\[
\int_{\Omega} \sigma (\mathbf{u}) : \varepsilon (\mathbf{v}-\mathbf{u}) \, dx = \int_{\Omega} \mathbf{f} \cdot (\mathbf{v}-\mathbf{u})\, dx + \int_{\Gamma} \sigma (\mathbf{u}) \mathbf{n} \cdot (\mathbf{v}-\mathbf{u})\, ds.
\]
From the definition of the Clarke subdifferential, the nonmonotone boundary condition (\ref{nbc}) is equivalent to 
\[
-\sigma_n (\mathbf{u})(v_n-u_n)\leq f^0(u_n;v_n-u_n).
\]
Here, the notation $f^0(x;z)$ stands for the generalized directional derivative of $f$ at $x$ in direction $z$. 
Substituting $\sigma(\mathbf{u})\mathbf{n}$ by $\mathbf{t}$ on $\Gamma_N$, using on $\Gamma_C$ the decomposition 
\[
\sigma (\mathbf{u})\mathbf{n} \cdot (\mathbf{v}-\mathbf{u}) = \sigma_t (\mathbf{u}) \cdot (\mathbf{v_t}-\mathbf{u_t}) +\sigma_n (\mathbf{u})(v_n-u_n)
\]
and taking into account  that on $\Gamma_C$  no tangential stresses are assumed, c.f.~\eqref{eq:notangstress},  
we obtain the  hemivariational inequality (HVI): Find $\mathbf {u} \in \mathcal{K}$ such that 
\begin{align} 
a(\mathbf{u},\mathbf{v}-\mathbf{u}) + \int_{\Gamma_C} f^0(u_n(s);v_n(s)-u_n(s)) \, ds \geq  \int_{\Omega} \mathbf{f} \cdot (\mathbf{v}-\mathbf{u})\, dx  \nonumber \\
+ \int_{\Gamma_N} \mathbf{t} \cdot (\mathbf{v}-\mathbf{u})\, ds \quad \forall \, \mathbf{v} \in \mathcal{K}. \label{HemVar}
\end{align}

\section{Boundary integral operator formulation}

Since the main difficulties of the boundary value problem, namely the nonmonotone adhesion law (\ref{nbc}) and the unilateral contact condition (\ref{unilateral}) appear on the boundary, 
it might be viewed advantageous to formulate \eqref{HemVar} as a boundary integral problem. 
To this end, we introduce the free boundary part $\Gamma_0=\Gamma \backslash \overline{\Gamma_D}=\Gamma_N \cup \Gamma_C $ and recall the Sobolev spaces \cite{HsWendl}: 
\[
\begin{array}{lll}
H^{1/2}(\Gamma)&=&\{v\in L^2(\Gamma)\, :\, \exists\, v' \in H^1(\Omega), \, \mbox{tr}\, v'=v\},\\[0.2cm]
H^{1/2}(\Gamma_0)&=&\{v=v'|_{\Gamma_0}  \, :\, \exists\, v'\in H^{1/2}(\Gamma)  \},\\[0.2cm]
\tilde{H}^{1/2}(\Gamma_0)&=&\{v=v'|_{\Gamma_0} \,:\, \exists\, v'\in H^{1/2}(\Gamma),\; \mbox{supp} \, v' \subset {\Gamma}_0 \} 
\end{array}
\]
with the standard norms
\[\|u\|_{H^{1/2}(\Gamma_0)}= \displaystyle\inf_{  v\in H^{1/2}(\Gamma), v|_{\Gamma_0}=u}\|v\|_{H^{1/2}(\Gamma)} \quad \text{and} \quad 
\|u\|_{\tilde{H}^{1/2}(\Gamma_0)}=\|u_0\|_{H^{1/2}(\Gamma)},
\]
where $u_0$ is the extension of $u$ onto $\Gamma$ by zero. The Sobolev space of negative order on $\Gamma_0$ are defined by duality as
\[
H^{-1/2}(\Gamma_0)=(\tilde{H}^{1/2}(\Gamma_0))^* \quad \text{and} \quad
\tilde{H}^{-1/2}(\Gamma_0) \; = \;  ({H}^{1/2}(\Gamma_0))^*.
\]
Moreover, from \cite[Lemma 4.3.1]{HsWendl} we have the inclusions
\[
\tilde{H}^{1/2}(\Gamma_0)\subset H^{1/2}(\Gamma_0)\subset L^2(\Gamma_0)\subset \tilde{H}^{-1/2}(\Gamma_0)\subset H^{-1/2}(\Gamma_0).
\]

For the solution $\mathbf{u}(\mathbf{x})$ of \eqref{eq1} with $\mathbf{x} \in \Omega\backslash \Gamma$ 
we have the following representation formula, also known as Somigliana's identity, see e.g.~\cite{Kleiber}:
{\small \begin{align} \label{eq2}
 \mathbf{u}(\mathbf{x})=\int_{\Gamma} G(\mathbf{x},\mathbf{y}) \left(\mathbf{T}_y \mathbf{u}(\mathbf{y})\right)\, ds_y -\int_\Gamma \mathbf{T}_y G(\mathbf{x},\mathbf{y})\mathbf{u}(\mathbf{y}) \,ds_y + \int_\Omega G(\mathbf{x},\mathbf{y})\mathbf{f}(\mathbf{y}) \, d\mathbf{y}, 
\end{align}
}
where $G(\mathbf{x},\mathbf{y})$ is the fundamental solution of the Navier-Lam\'{e} equation defined by
\begin{align*}
G(\mathbf{x},\mathbf{y})=\begin{cases}\frac{\lambda +3\mu}{4\pi \mu(\lambda+2\mu)}\left(\log |\mathbf{x}-\mathbf{y}|{\rm I} + \frac{\lambda + \mu}{\lambda +3\mu}\frac{(\mathbf{x}-\mathbf{y})(\mathbf{x}-\mathbf{y})^\top}{|\mathbf{x}-\mathbf{y}|^2}\right), \quad &\text{ if d=2}\\[0.2cm]
\frac{\lambda +3\mu}{8\pi \mu(\lambda+2\mu)}\left( |\mathbf{x}-\mathbf{y}|^{-1}{\rm I} + \frac{\lambda + \mu}{\lambda +3\mu}\frac{(\mathbf{x}-\mathbf{y})(\mathbf{x}-\mathbf{y})^\top}{|\mathbf{x}-\mathbf{y}|^3}\right), \quad &\text{ if d=3}
\end{cases}
\end{align*}
with the Lam\'{e} constants $\lambda, \mu >0$ depending on the material parameters, i.e.~the modulus of elasticity~$E$ and the Poisson's ratio~$\nu$:
$$
\lambda = \frac{E\nu}{1-\nu^2} \,, \qquad \mu= \frac{E}{1+\nu}\, . 
$$
 Here, $\mathbf{T}_y$  stands for the traction operator with respect to $\mathbf{y}$ defined by $\mathbf{T}_y (\mathbf{u}) := \sigma(\mathbf{u}(\mathbf{y}))\cdot \mathbf{n}_y$, where $\mathbf{n}_y$ is  the unit outer normal vector  at $\mathbf{y}\in \Gamma$. 
Letting $\Omega\backslash \Gamma \ni \mathbf{x} \to \Gamma$ in (\ref{eq2}), we obtain the well-known  
Calder\'{o}n operator 
\[
\left (\begin{array}{c}
          \mathbf{u}\\[0.2cm]
         \mathbf{T}_x \mathbf{u}
         \end{array}
\right )
= \left ( \begin{array} {cc}
           \frac{1}{2}I-K & V\\[0.2cm]
           W &      \frac{1}{2}I+K'  
           \end{array}
   \right) 
\left( \begin{array}{c}
\mathbf{u}\\[0.2cm]
\mathbf{T}_x \mathbf{u}  
\end{array}
\right ) + \left( \begin{array}{c} N_0 \mathbf{f} \\[0.2cm]
                                   N_1 \mathbf{f}
                   \end{array}
           \right),
\]       
with the single layer potential $V$, the double layer potential
$K$, its formal adjoint $K'$, and the hypersingular integral operator $W$ defined for $\mathbf{x} \in \Gamma$ as follows:
  \begin{alignat*}{2}
   \left( V \phi\right ) (\mathbf{x}) & :=  \int\limits_\Gamma G(\mathbf{x},\mathbf{y}) \phi (\mathbf{y}) \, ds_y, & \qquad
   \left(K \phi \right ) (\mathbf{x}) &:=
    \int_{\Gamma} \mathbf{T}_y G^\top(\mathbf{x},\mathbf{y})  \phi (\mathbf{y}) \, ds_y \\ 
   \left( K' \phi \right) (\mathbf{x}) & := \mathbf{T}_x \int\limits_\Gamma G 
    (\mathbf{x},\mathbf{y})  \phi (\mathbf{y}) \, ds_y, & \qquad
    \left(W \phi \right ) (\mathbf{x}) &:= - \mathbf{T}_x  \left (K \phi\right ) (\mathbf{x}).
  \end{alignat*}
  The Newton potentials $N_0, N_1$ are given for $x\in \Gamma$ by
\[
N_0 \mathbf{f} = \int\limits_\Gamma G 
    (\mathbf{x},\mathbf{y})  \mathbf{f} (\mathbf{y}) \, ds_y, \quad 
N_1 \mathbf{f} = \mathbf{T}_x \int\limits_\Gamma G 
    (\mathbf{x},\mathbf{y})  \mathbf{f} (\mathbf{y}) \, ds_y.
\]
From  \cite{Co-1} it is known that the linear operators
  \begin{alignat*}{2}
    V: &\;\mathbf{H} ^{-1/2 + \sigma} ( \Gamma)  \rightarrow  \mathbf{H}^{1/2 + \sigma} ( \Gamma)
    , &\qquad K:&\; \mathbf{H}^{1/2 + \sigma} ( \Gamma ) \rightarrow \mathbf{H}^{1/2 + \sigma} ( \Gamma) 
    \\
    K' : &\; \mathbf{H}^{-1/2 + \sigma} ( \Gamma )  \rightarrow  \mathbf{H}^{-1/2 + \sigma} (
    \Gamma) , &\qquad  
    W:&\; \mathbf{H}^{1/2 + \sigma} ( \Gamma ) \rightarrow \mathbf{H}^{-1/2 + \sigma} (\Gamma) 
  \end{alignat*}
are well-defined and continuous for $| \sigma | \leq \ds\frac{1}{2}$.   
Moreover, $V$ is symmetric and positive definite (elliptic on $\mathbf{H}^{-1/2}(\Gamma)$) in $\mathbb{R}^3$ and, if the capacity of $\Gamma$ is smaller than 1, also in $\mathbb{R}^2$. That can always be achieved by scaling, since the capacity (or conformal radius or transfinite diameter) of $\Gamma$ is smaller than 1, if $\Omega$ is 
contained in a disc with radius $<1$ 
(see e.g. \cite{SlSp,Steinbach}).  $W$ is symmetric and positive semidefinite with kernel $\R$ (elliptic on $\tilde{\mathbf{H}}^{1/2}(\Gamma_0)$). 
Hence, since $V : \mathbf{H}^{-1/2}(\Gamma) \to  \mathbf{H}^{1/2}(\Gamma)$ is invertible, we obtain by taking the Schur complement of the Calderon projector that
\begin{equation}\label{mainformula}
\mathbf{T}_x \mathbf{u}= P\mathbf{u} - N \mathbf{f}
\end{equation}
with the symmetric Poincar\'e-Steklov operator $P : \mathbf{H}^{1/2}(\Gamma) \to \mathbf{H}^{-1/2}(\Gamma)$ and its Newton potential $N$ given by
\begin{align*}
P\mathbf{u}  =  W \mathbf{u}  + \left( K' + \ds\frac{1}{2} I \right) V^{-1} \left( K + \ds\frac{1}{2} I \right) \mathbf{u} , \ 
N \mathbf{f} = \left( K' + \ds\frac{1}{2} I \right) V^{-1} N_0 \mathbf{f} -N_1 \mathbf{f}.
\end{align*}
Note that if $\mathbf{f}=0$, $P$ maps $\mathbf{u}$ to its traction and, therefore, the Poincar\'e-Steklov operator $P$ is sometimes called the Dirichlet-to-Neumann mapping. 
Moreover, $P$ induces a symmetric bilinear form on $\mathbf{H}^{1/2} ( \Gamma )$, and is continuous and $\tilde{\mathbf{H}}^{1/2}(\Gamma_0)$-elliptic, see e.g.~\cite{cf}, i.e.~there exist constants $c_P$, $C_P>0$ such that
\[
\begin{array}{rllll}
\|P \mathbf{v}\| _{\mathbf{H}^{-1/2}(\Gamma)}&\leq& C_P \|\mathbf{v}\| _{\mathbf{H}^{1/2}(\Gamma)} &\quad &\forall \,  \mathbf{v} \in {\mathbf{H}^{1/2}(\Gamma)}, \\[0.2cm]
\langle P\mathbf{v}, \mathbf{v} \rangle_{\Gamma_0} &\geq & c_P\|\mathbf{v}\|_{\mathbf{\tilde{H}}^{1/2}(\Gamma_0)} &\quad & \forall \,  \mathbf{v} \in {\mathbf{\tilde{H}}^{1/2}(\Gamma_0)}.
\end{array}
\]
Here, $\langle \cdot, \cdot \rangle_{\Gamma_0}$ is the $L^2(\Gamma_0)$-duality pairing between the involved spaces.  
Using the boundary functional sets 
$$
\mathcal{V}=\tilde{\mathbf{H}}^{1/2} (\Gamma_0) \quad \text{and} \quad \mathcal{K}^{\Gamma}=\{ \mathbf{v} \in \mathcal{V} : v_n \leq g \text{ a.e.~on } \Gamma_C \}
$$
we obtain as in the domain based cases the boundary integral hemivariational inequality (BIHVI), Problem ($\mathcal{P})$: Find $\mathbf{u} \in \mathcal{K}^{\Gamma}$ such that 
\begin{align}
\langle P \mathbf{u}, \mathbf{v}-\mathbf{u} \rangle_{\Gamma_0} + 
\int_{\Gamma_C} f^0(u_n(s);v_n(s)-u_n(s)) ds \geq  \langle \mathbf{t},\mathbf{v}-\mathbf{u} \rangle_{\Gamma_N} \nonumber \\ 
+ \langle N \mathbf{f}, \mathbf{v}-\mathbf{u} \rangle _{\Gamma_0} \quad \forall \,  \mathbf{v} \in \mathcal{K}^{\Gamma}. \label{BdHemVar}
\end{align}
To shorten the right hand side we introduce the linear functional 
\[
\langle \mathbf{F},\mathbf{v}\rangle_{\Gamma_0}=\int _{\Gamma_N} \mathbf{t} \cdot \mathbf{v}\, ds + \langle N \mathbf{f}, \mathbf{v} \rangle_{\Gamma_0}.
\]

\section{Regularization of the nonsmooth functional}

In this section, we review from \cite{Ovcharova2012,Ov_Gw} a class of smoothing approximations for nonsmooth functions that can be re-expressed by means of the  plus function $\mathrm{p}(t)=t^+=\max \{t, 0 \}$. The approximation is based   on  smoothing of the plus function via convolution. 


Let $\varepsilon >0$ be a small regularization parameter. 
The smoothing approximation $\hat{f}:\, \R \times \R_{++} \rightarrow \R$ of a locally Lipschitz function $f:\, \R \to \R$ is defined via convolution 
by
\[
\hat{f}(x, \varepsilon)=\int_{\R} f(x-\varepsilon t) \rho(t) \, dt,
\] 
where 
$\rho \, :\, \R \to \R_{+}$ is a probability density function such that
$$
\kappa= \int_{\R}|t|\, \rho(t)\, dt < \infty
$$
and
\[
\R_{+}=\{\varepsilon\in \R \, :\, \varepsilon \geq 0\}, \quad \R_{++}=\{\varepsilon\in \R \, :\, \varepsilon > 0\}.
\]
We consider a class of nonsmooth functions that can be reformulated by means of the plus function. 
To this class of functions belong   maximum, minimum or nested max-min function. 
 If $f$ is a maximum function, for example, $f(x)=\max\{g_1(x),g_2(x)\}$, then $f$ can be reformulated by using the plus function as 
 \begin{equation}\label{max2}
 f(x) = g_1(x) + \mathrm{p}( g_2(x)- g_1(x)),
 \end{equation}
 The smoothing function of $f$ is then defined by
 \begin{equation}\label{smooth_S}
S(x,\varepsilon) =g_1(x)+ \hat{p}(g_2(x)-g_1(x),\varepsilon), 
\end{equation}
where $\hat{p}( t, \varepsilon)$ is the smoothing function of $\mathrm{p}(t)$ via convolution.

  Choose, for example, the Zang probability density function
\[
\rho(t)= \left \{ \begin{array} {ll} 1 & \mbox{if} \, - \frac{1}{2}\leq t \leq
    \frac{1}{2} \\[0.2cm]
0 & \mbox{otherwise} \, .
\end{array} \right. 
\]
Then
\[
\hat{\mathrm{p}}(t,\varepsilon) = \int _\R \mathrm{p}(t-\varepsilon s)\rho(s) \, ds =\left \{ \begin{array} {ll} 0 & \mbox{if} \quad t < -\frac {\varepsilon}{2}\\ [0.1cm]
\frac{1}{2\varepsilon}(t+ \frac{\varepsilon}{2})^2 & \mbox{if} \, -
\frac{\varepsilon}{2} \leq t \leq \frac{\varepsilon}{2} \\[0.1cm]
t  &\mbox{if} \quad t > \frac
    {\varepsilon}{2}
\end{array}
\right. 
\]
and hence,
\begin{equation}\label{Zang}
{\small S(x,\varepsilon)= 
\left \{ \begin{array} {ll} g_1(x) & \; \mbox{if } g_2(x) -g_1(x) \leq - \frac{\varepsilon}{2} 
     \\[0.1cm]
\frac{1}{2\varepsilon} [g_2(x) -g_1(x)]^2  + \frac{1}{2} (g_2(x) +
g_1(x)) + \frac{\varepsilon}{8} & \; \mbox{if } |
g_2(x) -g_1(x)| \leq \frac{\varepsilon}{2} 
\\[0.1cm]
g_2(x) & \; \mbox{if } g_2(x) -g_1(x) \geq \frac{\varepsilon}{2} \; .
\end{array}
\right.}
\end{equation}
The relation (\ref{smooth_S}) 
can be extended to the maximum function $f:\R \to \R$ of $m$ continuous functions $g_1, \ldots, g_m$, i.e. 
\begin{align} \label{eq:bsp_f_function}
f(x)=\max \{g_1(x), g_2(x), \ldots, g_m(x)\}
\end{align}
by
\begin{equation} \label{S_general}
S(x, \varepsilon)= 
g_1(x) + \hat{\mathrm{p}}\left( g_2(x)-g_1(x)+\ldots + 
\hat{\mathrm{p}}\left(
g_m(x)-g_{m-1}(x), \varepsilon\right), \varepsilon\right).
\end{equation}
The major properties of the function $S$ in (\ref{S_general}) are listed in the following lemma:  
\begin{lemma}\cite[lemma 10] {Ovcharova2012}
Let $f$ and $S$ be defined as in \eqref{eq:bsp_f_function}, \eqref{S_general}, respectively, with $g_i \in C^1$. Then there holds:
\begin {description}
\item (i) For any $\varepsilon >0$ and for all $x\in \R$, 
\[
|S(x,\varepsilon)-f(x)|\leq (m-1)\kappa \varepsilon.
\]
\item(ii) The function $S$ is continuously differentiable on $\R\times \R_{++}$ and for any $x\in \R$ and $\varepsilon>0$ there exist $\Lambda_i \in [0,1]$ such that $\displaystyle \sum_{i=1}^m \Lambda_i=1$ and 
 \begin {equation}\label{f1}
  \frac{\partial S(x,\varepsilon)}{\partial x} = S_x(x, \varepsilon)=\displaystyle \sum_{i=1}^m \Lambda_i  g'_i(x).
\end{equation}
Moreover,
\begin{equation} \label{f1_1}
\left\{\lim_{z\to x, \varepsilon \to 0^+}  S_x(z, \varepsilon)\right\}\subseteq \partial f(x).
\end{equation} 
\end {description}
\end{lemma}
\begin{rem}
 Since $f(x)=\min\{g_1(x), \ldots, g_m(x)\} =  -\max \{-g_1(x), \ldots, -g_m(x)\}$ these type of functions can be handled as above.
\end{rem}
\begin{assumption}
Assume that there exists positive constants $c_i, d_i$ such that for all $x\in \R$
\begin{subequations}
 \begin{alignat}{2}
 | g'_i(x)| &\leq c_i(1+|x|) \label{ass_1}\\ 
 g'_i(x) \cdot (-x) &\leq d_i|x|. \label{ass_2}
\end{alignat}
\end{subequations}
\end{assumption}
By using (\ref{f1}) - (\ref{f1_1}) and (\ref{ass_1}) - (\ref{ass_2}), the following auxiliary result can be easily deduced. 
\begin{lemma} For any $x, z, \xi \in \R$ it holds that 
\begin{subequations}
\begin{alignat}{2}
\left| S_x (x, \varepsilon) \cdot z \right | &\leq  c(1+|x|)\,|z| \label{ass_3}\\
   S_x (x, \varepsilon) \cdot (-x)  &\leq d|x| \label{ass_4}\\
  \limsup_{z\to x, \varepsilon\to 0^+} S_x(\varepsilon,x)\, \xi &\leq f^0(x;\xi). 
\label{ass_5}
\end{alignat}
\end{subequations}
\end{lemma}
Further, we introduce  the functional $J_\varepsilon :\mathbf{H}^{1/2} (\Gamma) \to \R$ defined by
$$
J_\varepsilon(\mathbf{u})=\int_{\Gamma_C}  S( u_n(s), \varepsilon) \, ds.
$$
The regularized  problem $(\mathcal{P}_{\varepsilon})$ of ($\mathcal{P}$) is given by:~Find $\mathbf{u}_\varepsilon \in \mathcal{K}^{\Gamma}$ such that 
   \begin{equation} \label{reg} 
\langle P \mathbf{u}_{\varepsilon},\mathbf{v}-\mathbf{u}_{\varepsilon} \rangle_{\Gamma_0} + 
\langle D J_\varepsilon(\mathbf{u}_\varepsilon), \mathbf{v} -\mathbf{u}_\varepsilon \rangle_{\Gamma_C} 
 \geq  \langle \mathbf{F} , \mathbf{v}-\mathbf{u}_{\varepsilon} \rangle_{\Gamma_0} \quad \forall \,  \mathbf{v} \in \mathcal{K}^{\Gamma}, 
\end{equation} 
where $D J_\varepsilon :\mathbf{H}^{1/2} (\Gamma) \to \mathbf{H}^{1/2} (\Gamma)$ is the G\^{a}teaux derivative of the functional $J_\varepsilon$ and is given by
\[
\langle D J_\varepsilon(\mathbf{u}), \mathbf{v} \rangle_{\Gamma_C} = \int_{\Gamma_C} S_x ( u_n (s),\varepsilon) v_n(s)  \, ds.
\]
To simplify the notations we 
introduce  $\varphi: \mathbf{H}^{1/2} (\Gamma)\times \mathbf{H}^{1/2} (\Gamma) \to \R$ defined by 
\begin{equation} \label{varphi}
\varphi(\mathbf{u},\mathbf{v})=\displaystyle
{\int_{\Gamma_C}} f^0(u_n(s); v_n(s)- u_n(s)) \, ds \quad \forall \, 
\mathbf{u}, \mathbf{v} \in \mathbf{H}^{1/2} (\Gamma).
\end{equation}
The regularized version of $\varphi(\cdot,\cdot)$ is therewith $\varphi_{\varepsilon} : \mathbf{H}^{1/2} (\Gamma) \times \mathbf{H}^{1/2} (\Gamma) \to \R$, 
\[
\varphi_{\varepsilon}(\mathbf{u}, \mathbf{v}):=\langle D J_\varepsilon(\mathbf{u}), \mathbf{v}- \mathbf{u} \rangle_{\Gamma_C}. 
\]
The existence result for (\ref{BdHemVar}), resp. (\ref{reg}), is due to \cite{Gwinner_PhD} and relies in both cases on the pseudomonotonicity of $\varphi$ and  $\varphi_{\varepsilon}$, respectively. For details we refer the reader to \cite{Ovcharova2012,Ov_Gw}. 

Finally, we recall that the functional $\varphi : X\times X \to \R$, where $X$ is a real reflexive Banach space, is pseudomonotone if $u_n\rightharpoonup u$ (weakly) in $X$ and $\displaystyle \liminf _{n \to   \infty}  \varphi (u_n,u) \geq 0$ imply that, for all $v\in X$, we have $\limsup _{n \to \infty} \varphi (u_n,v)\leq \varphi (u,v).$

\section{Uniqueness results} \label{section4}
In this section, we discuss the uniqueness of solution of the boundary  hemivariational inequality 
and of the corresponding regularization problem. The main results are presented in Theorem \ref{theo2} which is based on the abstract uniqueness Theorem \ref{theo0} from \cite{Ovcharova2015} that  gives a sufficient condition for uniqueness. 
Similar uniqueness result but for the regularized problem is derived in Theorem \ref{theouniq}. We are also dealing  with the question which classes of smoothing functions preserve the property of unique solvability of the original problem as $\varepsilon \to 0^+$. 
\subsection{Uniqueness of the BIHVI}
\begin{assumption}
Assume that there exists an $\alpha_0 \in [0,c_P)$ such that for any $\mathbf{u}, \mathbf{v} \in \mathcal{V}$ it holds
\begin{equation} 
  \varphi(\mathbf{u},\mathbf{v})+\varphi(\mathbf{v},\mathbf{u})\leq \alpha_0 \|\mathbf{u}-\mathbf{v}\|^2_{\mathcal{V}}. \label{unique}
\end{equation}
To make the results of uniqueness self-consistent, we introduce from \cite{Ovcharova2015} the following theorem. 
\end{assumption}
\begin{theorem}\cite[Theorem 5.1]{Ovcharova2015}\label{theo0}
Under the assumption (\ref{unique}) with $\alpha_0<c_P$, there exists a unique solution to the BIHVI problem $(\mathcal{P})$, which depends Lipschitz continuously on the right hand side $\mathbf{F} \in \mathcal{V}^*$.
\end{theorem}
{\bf Proof}\, Assume that  $\mathbf{u}$, $\tilde{\mathbf{u}}$ are two solutions of $(\mathcal{P})$. Then the inequalities below hold:
\begin{align*}
\langle P\mathbf{u}-\mathbf{F}, \mathbf{v} -\mathbf{u}\rangle_{\Gamma_0} +
\varphi(\mathbf{u},\mathbf{v}) &\geq 0
\quad \forall \,  \mathbf{v}\in \mathcal{K}^{\Gamma} \\
\langle P \tilde{\mathbf{u}}-\mathbf{F}, \mathbf{v} -\tilde{\mathbf{u}}\rangle_{\Gamma_0} +
\varphi(\tilde{\mathbf{u}}, \mathbf{v}) &\geq 0
\quad \forall \,  \mathbf{v}\in \mathcal{K}^{\Gamma}.
\end{align*}
Setting $\mathbf{v}=\tilde{\mathbf{u}}$ in the first inequality and
$\mathbf{v}=\mathbf{u}$ in the second one, and summing up the resulting inequalities, we get
\begin{eqnarray} \label{unique01}
\langle P\mathbf{u} -P \tilde{\mathbf{u}}, \tilde{\mathbf{u}}
-\mathbf{u} \rangle_{\Gamma_0} +
\varphi (\mathbf{u},\tilde{\mathbf{u}}) +
\varphi (\tilde{\mathbf{u}},\mathbf{u}) \geq 0.
\end{eqnarray}
From the coercivity of the operator $P$ and the assumption (\ref{unique}) we obtain
$$
c_P\|\mathbf{u}-\tilde{\mathbf{u}}\|_{\mathcal{V}}^2 \leq \varphi (\mathbf{u} ,\tilde{\mathbf{u}}) +
\varphi (\tilde{\mathbf{u}} ,\mathbf{u} )\leq
  \alpha_0  \|\mathbf{u}-\tilde{\mathbf{u}}\|^2_{\mathcal{V}}.
$$
Hence, since $\alpha_0 \in [0, c_P)$,  if $\mathbf{u} \neq \tilde{\mathbf{u}}$ we receive a
contradiction. \\
Now let $\mathbf{F}_i \in \mathcal{V}^*$ and denote $\mathbf{u}_i=\mathbf{u}_{\mathbf{F}_i}, \; i=1,2.$  
Analogously to (\ref{unique01}), we find that
$$
\langle P\mathbf{u}_1 -\mathbf{F}_1-P \mathbf{u}_2+\mathbf{F}_2, \mathbf{u}_2
-\mathbf{u}_1  \rangle_{\Gamma_0}  +
\varphi (\mathbf{u}_1 ,\mathbf{u}_2) +
\varphi (\mathbf{u}_2 ,\mathbf{u}_1 ) \geq 0.
$$
Hence,
$$
c_P\|\mathbf{u}_1 -\mathbf{u}_2 \|^2_{\mathcal{V}}\leq \varphi (\mathbf{u}_1 ,\mathbf{u}_2 ) +
\varphi (\mathbf{u}_2 ,\mathbf{u}_1) + \langle \mathbf{F}_1-\mathbf{F}_2, \mathbf{u}_2-\mathbf{u}_1 \rangle_{\Gamma_0}
$$
and by (\ref{unique}), 
$$
(c_P-\alpha_0)\|\mathbf{u}_1-\mathbf{u}_2\|^2_{\mathcal{V}}\leq  \langle \mathbf{F}_1-\mathbf{F}_2, \mathbf{u}_2-\mathbf{u}_1 \rangle_{\Gamma_0} \leq \|\mathbf{F}_1-\mathbf{F}_2\|_{\mathcal{V}^*} \|\mathbf{u}_1 -\mathbf{u}_2 \|_{\mathcal{V}} \, .
$$
Also, since $\alpha_0 <c_P$ we deduce that
$$
\|\mathbf{u}_1 -\mathbf{u}_2 \|_{\mathcal{V}} \leq \frac{1}{c_P-\alpha_0} \, \|\mathbf{F}_1-\mathbf{F}_2\|_{\mathcal{V}^*},
$$
which concludes the proof of the theorem.
\qed

Further, we present  a class of locally Lipschitz functions for which the crucial assumption (\ref{unique}) is satisfied.
Let $f:\R \to \R$ be  a function such that
\begin{equation} \label{clmont}
( \xi^*-\eta^*)\,  (\xi- \eta)  \geq -\alpha_0 |\xi -\eta |^2 \quad 
\forall \,  \xi^*\in \partial f(\xi), \; \forall \,  \eta^*\in \partial f(\eta)
\end{equation}
for any $\xi, \eta \in \R$ and  some  $\alpha_0 \geq 0$.
From the definition of the Clarke generalized derivative \cite{Clarke} we get
$$
f^0(\xi;\eta-\xi)= \max _{\xi^*\in \partial f(\xi)}  \xi^* \,  (\eta -\xi). 
$$
Rewriting (\ref{clmont}) as

$$
\xi^* \, (\eta- \xi)  +  \eta^* \,  (\xi- \eta) \leq \alpha_0 |\xi -\eta |^2
$$
we find 
\begin{equation} \label{unique_f}
f^0(\xi;\eta-\xi)+ f^0(\eta;\xi-\eta) \leq \alpha_0|\xi-\eta|^2.
\end{equation}
Hence,
\begin{align}
 \varphi(\mathbf{u},\mathbf{v})+\varphi(\mathbf{v},\mathbf{u})&= 
\int_{\Gamma_C} f^0( u_n;  v_n - u_n )\, ds + \int_{\Gamma_C} f^0( v_n;  u_n - v_n )\, ds 
\nonumber \\
& \leq \alpha_0 \| u_n -  v_n\|^2_{L^2(\Gamma_C)} \leq \alpha_0  \|\mathbf{u}-\mathbf{v}\|^2_{\mathcal{V}} \label{est1_ex1}
\end{align}
and consequently the assumption (\ref{unique}) is satisfied provided that $\alpha_0$ is sufficiently small ($\alpha_0<c_P$).

Next, we show that if $\partial f$ includes only non-negative jumps, then the condition (\ref{clmont}) is globally satisfied. Whereas for negative jumps, (\ref{clmont}) holds only locally. 
\begin{figure}[t] 
\centering
\includegraphics[trim=0cm 0cm 0cm 0cm,clip,width=0.6\textwidth]{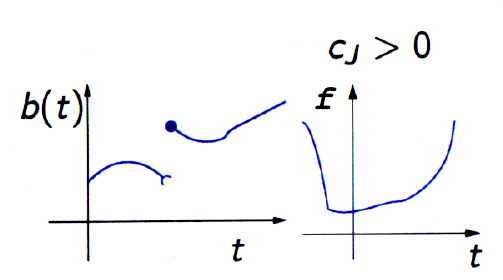} 
\caption{An exemplary function $b$ with one positive jump and its anti-derivative  $f(x)= \int_0^x b(t)\,dt$ }
 \label{nonnegative_jump}
\end{figure}
\begin{example} \label{ex1}
Let $f :\R \to \R$ be a  continuous, piecewise $C^{0,1}$ function such that its first derivative $b(x):=f'(x)$ has finite non-negative jumps at the points $x_i^J$ of discontinuity (see Fig. \ref{nonnegative_jump}). This means that 
$\underline{b}(x_i^J)\leq \overline{b}(x_i^J)$, where  
$$
\underline{b}(x_i^J):=f'(x_i^J-0)= \lim _{h\to 0^-}\frac{f(x_i^J+h)-f(x_i^J)}{h}
$$
and 
$$
\overline{b}(x_i^J):= f'(x_i^J+0)=\lim _{h\to 0^+}\frac{f(x_i^J+h)-f(x_i^J)}{h}.
$$
Hence, 
$$
\partial f(x_i^J) =\left[\underline{b}(x_i^J), \, \overline{b}(x_i^J) \right].
$$
\\
Setting $\mathcal{E}=\cup\,  \mathcal{E}_i$ the finite set of open intervals $\mathcal{E}_i$, where  the function $b|_{\mathcal{E}_i}$ is smooth with the Lipschitz constant $c_{L_i}$, we define
\begin{equation} \label{eq1_ex1}
\alpha_0=\displaystyle \max_{i}  \, c_{L_i} . 
\end{equation} 
Let $\{x_i^J\}_{i=1}^k$  be the set of jags between $x_1$ and $x_2$  with $ x_2 < x_1^J < \cdots < x_k^J<x_1 $, where $\underline{b}(x_i^J)\leq \overline{b}(x_i^J)$.  Then, for any $x_1>x_2$ we have the following
\begin{eqnarray*}
\frac{b(x_1)-b(x_2)}{x_1-x_2} & = & \frac{b(x_1)-b(x_2)+ \displaystyle \sum _{i=1}^k\overline{b}(x^J_i)-\displaystyle \sum _{i=1}^k\overline{b}(x^J_i)}{x_1-x_2} \\
&\geq& \frac{b(x_1)-b(x_2)+ \displaystyle \sum _{i=1}^k\underline{b}(x^J_i)-\displaystyle \sum _{i=1}^k\overline{b}(x^J_i)}{x_1-x_2}\\
&= &\frac{b(x_1)-\overline{b}(x^J_k) -b(x_2)+\underline{b}(x^J_1)
+ \displaystyle \sum _{i=2}^k \big(\underline{b}(x^J_i)-\overline{b}(x^J_{i-1})\big)}{x_1-x_2} \\
&\geq& \frac{-\alpha_0(x_1-x^J_k)-\alpha_0(x^J_1-x_2)-\alpha_0 \displaystyle \sum _{i=2}^k  (x^J_i -x^J_{i-1})}{x_1-x_2}\\ 
&=& \frac{-\alpha_0(x_1-x_2)}{x_1-x_2}=-\alpha_0,
\end{eqnarray*}
from which the assumption (\ref{clmont}) follows immediately.

\qed 
\end{example}
\begin{rem}
If the graph of  $\partial f$ consists  of several decreasing straight line  segments and non-negative jumps,  
then the value $-\alpha_0$ in (\ref{clmont}) is the steepest decreasing slope.
\end{rem}
The next example treats the case of negative jumps.
\begin{example} \label{ex2}
Let $f :\R \to \R$ be a  continuous, piecewise $C^{0,1}$ function such that its first derivative  has finite negative jumps at the points $x_i^J$ of discontinuity (see Fig. \ref{negative_jump}). 
In this case, 
$$
\partial f(x_i^J) =\left[\overline{b}(x_i^J), \underline{b}(x_i^J) \right].
$$
Then, for any $x_1>x_2$, we have
\begin{equation} \label{ex2_0}
\frac{b(x_1)-b(x_2)}{x_1-x_2} \geq -\alpha_0 + \frac{\displaystyle \sum _{i=1}^k c^J_i}{x_1-x_2},
\end{equation}
with $\alpha_0\geq 0$ as defined in (\ref{eq1_ex1}) 
and $c^J_i=\overline{b}(x^J_i)-\underline{b}(x^J_i)<0$ standing for the value of the jump in the point $x^J_i$. \\
\begin{figure}[t] 
\centering
\includegraphics[trim=0cm 0cm 0cm 0cm,clip,width=0.6\textwidth]{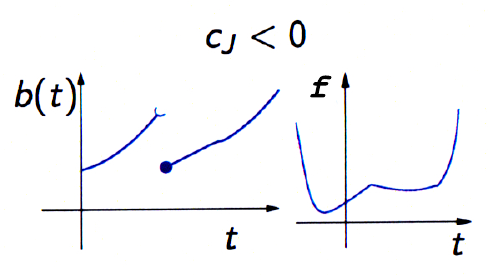} 
\caption{An exemplary function $b$ with one negative jump and its anti-derivative $f(x)=\int_0^x b(t) \, dt$}
 \label{negative_jump}
\end{figure}
Further, for the sake of simplicity, we consider a function with one negative jump $c_J$  and estimate as follows:
\begin{equation} \label{ex2_1}
c_J(x_1-x_2) =   c_J(x_1-x_J) + c_J(x_J-x_2)  =  c_J(x_1-x_J)_{+} + c_J(x_J-x_2)_{+}. 
\end{equation}
Using $x^2\geq (x-a)_{+}$ for any $a\geq \frac{1}{4}$, we get
$$
(x_1-x_2)^2\geq (x_1-x_2-\frac{1}{4})_{+}\geq (x_1-x_2-(x_J-x_2))_{+}=(x_1-x_J)_{+} \quad \forall \,  x_2\leq x_J -\frac{1}{4}
$$
as well as
$$
(x_1-x_2)^2\geq (x_1-x_2-\frac{1}{4})_{+}\geq (x_1-x_2-(x_1-x_J))_{+}=(x_J-x_2)_{+} 
\quad \forall \,  x_1\geq x_J +\frac{1}{4}.
$$
Hence, the condition  (\ref{clmont}) is fulfilled with  $\alpha_0 -2c_J>0$, but only for those  $x_1, x_2$ satisfying $x_1\geq x_J + \frac{1}{4}$ and $x_2 \leq x_J - \frac{1}{4}$.

Let $h(t)$ be a Lipschitz continuous function with Lipschitz constant $c_L$. We  include the negative jump $c^J$ at the point $t_J$ through the Heaviside function $\Theta$ by 
\[
b(t)=h(t)+c^J\Theta(t-t_J)
\]
and compute
 \[
 f^0(t;d)=\left \{ \begin{array} {ll}
                            h(t)d, & \, \mbox{if} \; t<t_J\\[0.1cm]
                            h(t)d+c_Jd, & \, \mbox{if} \; t> t_J,
                            \end{array}
                     \right.
                     \quad f^0(t_J;d)=\left \{ \begin{array} {ll}
                            h(t_J)d, & \, \mbox{if} \; d>0\\[0.1cm]
                            h(t_J)d+c_Jd, & \, \mbox{if} \; d\leq 0.
                            \end{array}
                     \right.
      \] 
Hence, for any $t_1, t_2 <t_J$, we have
\begin{align} 
f^0(t_1;t_2-t_1)&+ f^0(t_2;t_1-t_2)= h(t_1)(t_2-t_1)+h(t_2)(t_1-t_2) \nonumber \\[0.1cm]
=&(h(t_2)-h(t_1))(t_1-t_2)\leq  c_L(t_1-t_2)^2. \label{est1_ex2}
\end{align}
Analogously, for any $t_1, t_2 > t_J$, we get
\begin{align*} 
f^0(t_1;t_2-t_1)+f^0(t_2;t_1-t_2)=& h(t_1)(t_2-t_1)+c_J(t_2-t_1)+h(t_2)(t_1-t_2)\\[0.1cm]
&+c_J(t_1-t_2)\leq c_L(t_1-t_2)^2.
\end{align*}  
In both cases, the condition (\ref{unique_f}) is satisfied.     
\\
Let now $t_1< t_J \leq t_2$, then we have
\begin{align} \label{est2_ex2}
 f^0(t_1;t_2-t_1)&+f^0(t_2;t_1-t_2)= h(t_1)(t_2-t_1)+h(t_2)(t_1-t_2) + c_J(t_1-t_2) \nonumber \\[0.1cm]
 & \quad \leq c_L(t_1-t_2)^2 + c_J(t_1-t_2) \leq \alpha_0 (t_1-t_2)^2
\end{align} 
for some $\alpha_0$ if and only if
\begin{align} \label{est3_ex2}
\frac{c_J}{t_1-t_2}+c_L \leq \alpha_0.
\end{align}
 Thus, we assume 
\begin{align} \label{est4_ex2}
c_L<\alpha_0 \quad \mbox{and} \quad t_1 \leq t_J- \frac{c_J}{c_L-\alpha_0},
\end{align}
which, since $t_J \leq t_2$, implies (\ref{est3_ex2}) immediately. \\

Based on the results presented in the Examples~\ref{ex1} and \ref{ex2} we derive the following uniqueness result.  
\begin{theorem} \label{theo2}
Let $\mathbf{u}$ be a solution of the BIHVI problem $(\mathcal{P})$. Then $\mathbf{u}$ is unique if one of the following conditions holds:
\begin{itemize}
\item [(a)] The jumps $c^i_J$ are non-negative and 
the Lipschitz constant $c_L$ satisfies $c_L < c_P$.
\item [(b)] The jump $c_J$ is negative, $c_L<c_P$  
and  ${u}_n(x)\leq t_J-\frac{1}{\tilde{c}}$ on $\Gamma_C$ 
for some positive $\tilde{c}$ such that $\tilde{c}<-\frac{c_P-c_L}{c_J}$. 
\end{itemize}
\end{theorem}
{\bf Proof} (b)\, Let $\mathbf{u} \in \mathcal{V}$ be such that ${u}_n(x)\leq t_J-\frac{1}{\tilde{c}}$ on $\Gamma_C$  and  $\mathbf{v} \in \mathcal{V}$ arbitrary. With $\alpha_0:= c_L-\tilde{c}\,c_J$, using  (\ref{est1_ex2}) and (\ref{est2_ex2}), we obtain
\begin{align}
\varphi(\mathbf{u}, \mathbf{v}) & + \varphi(\mathbf{v}, \mathbf{u}) =  \int_{\Gamma_C } f^0(u_n;v_n-u_n) \, ds + \int_{  \Gamma_C } f^0(v_n;u_n-v_n) \, ds \nonumber \\
=& \int_{ \{ s\in \Gamma_C \,:\, v_n(s)<t_J \}} f^0(u_n;v_n-u_n) \, ds + \int_{ \{ s\in \Gamma_C \,:\, v_n(s)\geq t_J \}} f^0(u_n;v_n-u_n) \, ds \nonumber \\ 
& +  \int_{ \{ s\in \Gamma_C \,:\, v_n(s)<t_J \}} f^0(v_n;u_n-v_n) \, ds  + \int_{ \{ s\in \Gamma_C \,:\, v_n(s)\geq t_J \}} f^0(v_n;u_n-v_n) \, ds \nonumber \\
  \leq & \int_{ \{ s\in \Gamma_C \,:\, v_n(s)<t_J \}} c_L (u_n-v_n)^2 \, ds + \int_{ \{ s\in \Gamma_C \,:\, v_n(s) \geq t_J \}} \alpha_0 (u_n-v_n)^2 \, ds \nonumber \\
  \leq & \max \{c_L, \alpha_0\} \left\|u_n-v_n \right\|^2_{L^2(\Gamma_C)} \leq \alpha_0 \left\|\mathbf{u}-\mathbf{v} \right\|^2_{\mathcal{V}}.  \label{est5_ex2}
\end{align}
Further, we assume that $\mathbf{\tilde{u}} \in \mathcal{V}$ is another solution of $(\mathcal{P})$. Following the proof of Theorem \ref{theo0}, in virtue of  $\alpha_0 = c_L-\tilde{c}\,c_J < c_p$ and (\ref{est5_ex2}) with $\mathbf{v}:=\mathbf{\tilde{u}}$, we obtain a contradiction. \qed
\begin{rem}
We note that our uniqueness result is more general then the result in \cite[Theorem 7]{Nesemann}, which demands $\tilde{c}=4$.  
\end{rem}

\end{example}
 
\subsection{Uniqueness of the regularized problem}
\begin{assumption}
Assume that the assumption (\ref{clmont}) is satisfied for the regularization function, i.e.~there exists a constant $\alpha_0 \geq 0$  (in general depending on $\varepsilon>0$) such that
\begin{equation} \label{uniqreg}
\left (  S_x(x_1, \varepsilon)-  S_x(x_2, \varepsilon) \right )  (x_1-x_2) \geq - \alpha_0|x_1-x_2|^2 \quad \forall \,  x_1, x_2 \in \R.
\end{equation}
\end{assumption}
Hence, for any $\mathbf{u}, \mathbf{v} \in \mathcal{V}$, we have
\begin{align} 
 \langle D J_{\varepsilon}(\mathbf{u}) - D J_{\varepsilon}(\mathbf{v}),\mathbf{v}-\mathbf{u}\rangle_{\Gamma_C} 
&=  \int_{\Gamma_C}\left (  S_x( {u}_n(s),\varepsilon) - S_x(
{v}_n(s),\varepsilon)\right )(v_n(s)- {u}_n(s))\,ds \nonumber \\
&\leq  \alpha_0 \| {u}_n-  {v}_n\|^2_{L^2(\Gamma_C)} \leq \alpha_0 \|\mathbf{u}-\mathbf{v}\|^2_{\mathcal{V}}. \label{ass1_th5}
\end{align}
Due to  Theorem~\ref{theo0}, we have the following uniqueness result for the regularized problem. 
\begin{theorem}\label{theouniq}
Under the assumption \eqref{uniqreg} with $\alpha_0 < c_p$, there exists a unique solution to the regularized problem $(\mathcal{P}_{\varepsilon})$, which depends Lipschitz continuously on the right hand side $\mathbf{F} \in \mathcal{V}^*$.
\end{theorem}
One of the most interesting question is, which classes of smoothing functions preserve the property of unique solvability of the original problem as $\varepsilon \to 0^+$. As the next example shows, this is the case for the smoothing function (\ref{Zang}) relevant to our computation. In particular, we  show that the constant $\alpha_0$ in (\ref{uniqreg}) does not depend on $\varepsilon$. 
\begin{example} \label{ex3} In view of our application to contact problems with nonmonotone laws in the boundary conditions, we consider $f(x)=\max\{g_1(x), g_2(x)\}$  with 
$
g_i(x)=\frac{a_i}{2}x^2 +c_ix +d_i. 
$
We use the smoothing function (\ref{Zang}) as a regularization of $f$.  The corresponding approximation of $\partial f$  is given by
\[
 S_x(x, \varepsilon) =  \left \{ \begin{array} {ll}
 g_1'(x) & \mbox{if} \, \;  g_2(x) -g_1(x) \leq - \frac{\varepsilon}{2} \\[0.1cm] 
 \lambda_1  g_1'(x) +\lambda_2  g_2'(x) 
& \mbox{if} \, \;  \left |g_1(x) -g_2(x) \right | < \frac{\varepsilon}{2} \\[0.1cm] 
 g_2'(x) & \mbox{if} \, \; g_2(x) -g_1(x) \geq \frac{\varepsilon}{2},
\end{array}
\right.
\]
where
$$
\lambda_2= \frac{1}{\varepsilon} (g_2(x) -g_1(x))+ \frac{1}{2},\quad  \lambda_1 =
1 -  \lambda_2.
$$
We need to consider only the region $I_{\varepsilon}:=\{x\in \R \, :\, |g_1(x)-g_2(x)|<\frac{\varepsilon}{2}\}$, where $f$ is actually approximated. Otherwise, the condition    (\ref{uniqreg}) is automatically satisfied. \\ 
Let $x_i \in I_{\varepsilon}, \; i=1,2$. Setting $g(x):=g_1(x)-g_2(x)$, we introduce
$$
\lambda_2^i= \frac{1}{\varepsilon}(g_2(x_i)-g_1(x_i))=-\frac{1}{\varepsilon}g(x_i) \; 
\mbox{and} \; \lambda_1^i=1-\lambda_2^i,\; i=1,2. 
$$
Hence,
\begin{equation} \label {defg}
\lambda_1^1-\lambda_1^2 =  \lambda_2^2-\lambda_2^1 
= \frac{1}{\varepsilon} (g(x_1)-g(x_2))
\end{equation}
and thus,
\begin{align}
 S_x(x_1, \varepsilon) - S_x(x_2, \varepsilon) &=  \lambda_1^1(a_1 x_1+c_1)+  \lambda_2^1(a_2 x_1 +c_2) 
 -\lambda_1^2(a_1 x_2 +c_1)  \nonumber\\
 & \quad -  \lambda_2^2(a_2 x_2 +c_2) \nonumber \\ 
 & = a_1 \lambda_1^1 (x_1-x_2)+a_2 \lambda_2^1 (x_1-x_2) 
  + (\lambda_1^1-\lambda_1^2)(a_1x_2+c_1)   \nonumber\\
 & \quad  + (\lambda_2^1-\lambda_2^2)(a_2 x_2+c_2)\nonumber\\
 & = \left (  a_1 \lambda_1^1 + a_2 \lambda_2^1 \right ) (x_1-x_2) 
  +(\lambda_1^1-\lambda_1^2) (g_1^{\prime}(x_2)- g_2^{\prime}(x_2))\nonumber\\
 &= \left (a_1 \lambda_1^1+a_2 \lambda_2^1 \right )(x_1-x_2)+ (\lambda_1^1-\lambda_1^2) g^{\prime}(x_2). \label{onejag1}
\end{align}
Similarly to (\ref{onejag1}), we find that
\begin{equation} \label{onejag2}
\hspace{-1cm}S_x(x_1, \varepsilon) - S_x(x_2, \varepsilon)= \left( a_1 \lambda_1^2 + a_2 \lambda_2^2\right) (x_1-x_2) 
+ (\lambda_1^1-\lambda_1^2) g^{\prime}(x_1).
\end{equation}
Summing now (\ref{onejag1}) and (\ref{onejag2}),  we get 
\begin{align} \label{onejag3}
2\left (S_x(x_1, \varepsilon) - S_x(x_2, \varepsilon)\right ) &= (a_1 \lambda_1^1+a_2 \lambda_2^1+ a_1 \lambda_1^2 + a_2 \lambda_2^2) (x_1-x_2)  \nonumber\\
 & \quad 
+ (\lambda_1^1-\lambda_1^2) ( g^{\prime}(x_1) +g^{\prime}(x_2) ). 
\end{align}
From (\ref{defg}), using the structure of $g_1$ and $g_2$, we find that 
\begin{align}
\lambda_1^1-\lambda_1^2 &=  \frac{1}{\varepsilon} (g(x_1)-g(x_2)) =  \frac{1}{\varepsilon} \left( \frac{a_1-a_2}{2}(x_1^2-x_2^2) +(c_1-c_2)(x_1-x_2)\right )\nonumber \\ 
&=  \frac{1}{\varepsilon} (x_1-x_2) \left( \frac{a_1-a_2}{2} x_1 + \frac{c_1-c_2}{2}+ \frac{a_1-a_2}{2} x_2 + \frac{c_1-c_2}{2} \right )\nonumber \\ 
&=  \frac{1}{2\varepsilon} (x_1-x_2) ( g^{\prime}(x_1) +g^{\prime}(x_2) ). \label{onejag4}
\end{align}
Hence,
\begin{align*}
 & 2\left (S_x(x_1, \varepsilon)  - S_x(x_2, \varepsilon)\right ) (x_1-x_2) 
 \geq (a_1 \lambda_1^1+a_2 \lambda_2^1+ a_1 \lambda_1^2 + a_2 \lambda_2^2) (x_1-x_2)^2 \\
 & \quad =(a_1 (1-\lambda_2^1)+a_2 \lambda_2^1+ a_1 (1-\lambda_2^2) + a_2 \lambda_2^2) (x_1-x_2)^2\\
 & \quad \geq (-|a_1| |(1-\lambda_2^1)| - |a_2| |\lambda_2^1|- |a_1| |(1-\lambda_2^2)| - |a_2| |\lambda_2^2|) (x_1-x_2)^2\\
 & \quad \geq -\max\{|a_1| , |a_2|\} ( |(1-\lambda_2^1)| + |\lambda_2^1|+  |(1-\lambda_2^2)| + |\lambda_2^2|) (x_1-x_2)^2\\
 & \quad \geq -\max\{|a_1| , |a_2|\} ( 3/2 + 1/2+  3/2 + 1/2) (x_1-x_2)^2\\
 & \quad =-4\max\{|a_1| , |a_2|\} (x_1-x_2)^2
\end{align*}
and consequently, (\ref{uniqreg}) is satisfied with $\alpha_0 = 2\max\{|a_1| , |a_2|\}$.  
\end{example}
\begin{rem}
If we use  (\ref{Zang}) as a smoothing approximation for the maximum (resp. minimum) function  of several quadratic functions   $a_i x^2 +c_i x +d_i$, then (\ref{uniqreg}) is satisfied and, therefore, the regularized problem is always uniquely solvable if all $|a_i|$ are sufficiently small. Since the constant $\alpha_0$ does not depend on $\varepsilon$,  the property of uniqueness is preserved as $\varepsilon \to 0^+$.  
\end{rem}

\section{Discretization with boundary elements}

To avoid domain approximation, let $\Omega\subset \R^{d}$, $d=2,3$,  also be a polygonal domain. Let $\mathcal{T}_h$ be a sufficiently fine finite element mesh of the boundary~$\Gamma$ respecting the decomposition of $\Gamma$ into $\Gamma_D$, $\Gamma_N$ and $\Gamma_C$, 
$p=(p_T)_{T\in \mathcal{T}_h}$ a polynomial degree distribution over $\mathcal{T}_h$, $\mathbb{P}_{p_{T}}(\hat{T})$ the space of polynomials of order $p_{T}$ on the reference element $\hat{T}$, and $\Psi_T: \hat{T} \rightarrow T \in \mathcal{T}_h$ a bijective, (bi)-linear transformation. In 2D, $\hat{T}$ is the interval $[-1,1]$, whereas in 3D it is the reference square $[-1,1]^2$. Let $\Sigma_{T, hp}$ be the set of all affinely transformed (tensor product based) Gauss-Lobatto nodes on the element $T$ of the partition $\mathcal{T}_h$ of $\Gamma$, corresponding to the polynomial degree $p_T$, and set $\displaystyle \Sigma_{hp}:=\displaystyle \bigcup_{T\in \mathcal{T}_h|_{\Gamma_C}}\Sigma_{T, hp}$, see \cite{krebs2007p,MS-1,gwinner2013hp}. Furthermore, we assume in this section that $g \in C^0(\overline{\Gamma_C})$ to allow point evaluation.\\
For the discretization of the displacement $\mathbf{u}$ we use 
\begin{align*}
\mathcal{V}_{hp}&=\{\mathbf{v_{hp}} \in \mathbf{C}^0(\Gamma)\, :\, \mathbf{v_{hp}}|_{T} \circ \Psi_T  \in [\mathbb{P}_{p_{T}
}(\hat{T})]^{d} \quad \forall \,  T \in \mathcal{T}_h, \;  \mathbf{v_{hp}} =0 \; \mbox{on} \;  \overline{\Gamma_D} \} , \\
\mathcal{K}_{hp}^{\Gamma} &=\{ \mathbf{v_{hp}} \in \mathcal{V}_{hp} \, :\, (\mathbf{v_{hp}}\cdot \mathbf{n})(P_i) \leq g(P_i) \quad \forall \,  P_i \in \Sigma_{hp} \}.
\end{align*}
In general $\mathcal{K}_{hp}^{\Gamma} \not\subseteq   \mathcal{K} ^{\Gamma}$. For the approximation of the Poincar\'e-Steklov operator, namely $V^{-1}$, we need the space 
\[
\mathcal{W}_{hp}= \{ \psi_{hp} \in \mathbf{L}^2(\Gamma) \, : \, \psi_{hp}|_{T} \circ \Psi_T \in [\mathbb{P}_{p_T-1}(\hat{T})]^{d} \quad \forall \,  T \in \mathcal{T}_h\} \subset \mathbf{H}^{-1/2}(\Gamma).
\]
For more details on the approximation techniques based on boundary element method see e.g. \cite{ccjg,Co-1,CoSt,Gu,GS-1,Ha,MS-1,MS-2,SlSp}. 
Let $\{\phi_i\}_{i=1}^{N_D}$ and $\{\psi_j\}_{j=1}^{N_N}$ be a bases of $\mathcal{V}_{hp}$ and $\mathcal{W}_{hp}$, respectively. Then, the  boundary matrices are given by
\begin{align*}
(V_{hp})_{j,i}  &= \langle V \psi_i, \psi_j\rangle, \quad  (K_{hp})_{j,i} = \langle K \phi_i, \psi_j\rangle, \\
(W_{hp})_{j,i} &= \langle W \phi_i, \phi_j\rangle, \quad
(I_{hp})_{j,i}  = \langle \phi_i, \psi_j\rangle 
\end{align*}
and, therewith, the approximation of the Galerkin matrix is
\[
{P}_{hp}=W_{hp}+\left( K_{hp}+ \frac{1}{2}I_{hp}\right)^\top V^{-1}_{hp} \left( K_{hp}+ \frac{1}{2}I_{hp}\right).
\] 
With the canonical embeddings
\begin{gather*}
 i_{hp} :  \mathcal{W}_{hp}  \hookrightarrow \mathbf{H}^{-1/2}(\Gamma) \, , \qquad
 j_{hp} :  \mathcal{V}_{hp}  \hookrightarrow \mathbf{H}^{1/2}(\Gamma)
\end{gather*}
and their duals $i_{hp}^*$ and $j_{hp}^*$, 
the discrete Poincar\'{e}-Steklov operator $P_{hp} : \mathcal{V}_{hp} \to \mathcal{V}_{hp}^*$ can be  also represented   by 
\[
P_{hp}= j_{hp}^*W j_{hp}+ j_{hp}^*\left(K'+\frac{1}{2}I\right )i_{hp} (i_{hp}^* V i_{hp})^{-1} i_{hp}^*\left(K+\frac{1}{2}I\right)j_{hp}. 
\]
According to \cite{cc}, see e.g.~\cite{Banz2015Stab} for the $hp$-version, this functional is well-defined and satisfies
\begin{equation} \label{coerc}
\exists\, c >0: \quad \langle P_{hp} \mathbf{v}_{hp}, \mathbf{v}_{hp} \rangle_{\Gamma_0} \geq c \| \mathbf{v}_{hp}\|_{\tilde{\mathbf{H}}^{1/2}(\Gamma_0)} \qquad \forall\,  \mathbf{v}_{hp} \in \mathcal{V}_{hp}
\end{equation}
if $\min\{h,p^{-1}\}$ is sufficiently small.
Further, we consider the operator 
$E_{hp}: \mathbf{H}^{1/2}(\Gamma)\to \mathbf{H}^{-1/2}(\Gamma)$, reflecting the consistency error in the discretization of the Poincar\'{e}-Steklov operator $P$, defined by
\[
E_{hp}:=P-P_{hp}=\left(\frac{1}{2}I+K' \right)\left(V^{-1}-i_{hp}(i_{hp}^*Vi_{hp})^{-1}i_{hp}^*\right)\left(\frac{1}{2}I+K\right). \]
From \cite{MS-1} (see also \cite{cc}), the operator $E_{hp}$ is bounded and there exists a constant $c>0$ such that
\[
\|E_{hp}\mathbf{v}\|_{\mathbf{H}^{-1/2}(\Gamma)}\leq c \, \displaystyle \inf_{{\psi_{hp}}\in \mathcal{W}_{hp}}\|V^{-1}(I+K)\mathbf{v} - {\psi_{hp}}\|_{\mathbf{H}^{-1/2}(\Gamma)} \quad \forall \,  \mathbf{v} \in \mathbf{H}^{1/2}(\Gamma).
\]
Now, we turn to the discretization of the regularized problem $(\mathcal{P}_\varepsilon)$. 
The discretized regularized  problem ($\mathcal{P}_{\varepsilon, hp}$) is: Find
   $\mathbf{u}^{\varepsilon}_{hp} \in \mathcal{K}^{\Gamma}_{hp}$ such that for all $\forall \,  \mathbf{v}_{hp} \in \mathcal{K}_{hp}^{\Gamma}$ 
\begin{equation} \label{reg_dis_problem}
\langle {P}_{hp} \mathbf{u}^{\varepsilon}_{hp}, \mathbf{v}_{hp}-\mathbf{u}^{\varepsilon}_{hp} \rangle_{\Gamma_0}  +
\langle D J_{\varepsilon}(\mathbf{u}^{\varepsilon}_{hp}), \mathbf{v}_{hp} -\mathbf{u}^{\varepsilon}_{hp} \rangle_{\Gamma_C} 
 \geq \langle \mathbf{F}, \mathbf{v}_{hp}-\mathbf{u}^{\varepsilon}_{hp} \rangle_{\Gamma_0}.   
\end{equation}
The convergence result for the $hp$-solution $\mathbf{u}^\varepsilon_{hp}$ of ($\mathcal{P}_{\varepsilon, hp}$) is due to the following abstract approximation result from \cite{Gw_Ov}. 

Let $\mathcal{K}$ be  a closed, convex nonempty subset of a reflexive Banach space $X$ and  $VI(\psi,f,\mathcal{K})$ be the pseudo-monotone variational inequality:~Find $u \in \mathcal{K}$ such that 
$
\psi(u,v) \geq \langle f, v-u\rangle \, \forall \,  v \in \mathcal{K}.
$
Let $T$ be a directed set.  We  introduce the family $\{\mathcal{K}_t\}_{t\in T}$ of nonempty, closed and convex sets and admit the  following hypotheses:
\begin{description}
\item (H1) If  $\{v_{t'}\}_{t'\in T'}$ weakly converges to $v$ in $X$, $v_{t'} \in
  \mathcal{K}_{t'} ~(t' \in T') $ for a subnet $\{\mathcal{K}_{t'}\}_{t'\in T'}$ of the  net  $\{\mathcal{K}_t\}_{t\in T}$,  then   $v\in \mathcal{K}$.
\item (H2) For any $v\in \mathcal{K}$ and any $t \in T$ there exists $v_{t} \in \mathcal{K}_{t}$ 
such that  $v_{t} \to v$  in $X$.
\item (H3) $\psi_t$ is pseudomonotone for any $t \in T$.
\item (H4) $f_t \to f$ in $V^\ast$.
\item (H5) For any nets $\{u_t\}$ and $\{v_t\}$ 
 such that  $u_t \in \mathcal{K}_t$,   $v_t\in \mathcal{K}_t$, 
  $u_t \rightharpoonup u$, and $v_t \to v$ in $X$ it follows that 
$$
\displaystyle {\liminf_{t \in T}}\,~ \psi_t(u_t,v_t) \leq \psi(u,v) \,.
$$
\item  (H6) 
  There exist constants $c>0$,   $d$, $d_0 \in \R$ and $\alpha >1$  (independent of   $t \in T$) such that for some $w_t \in \mathcal{K}_t$
with $w_t \to w$ there holds  
$$
-\psi_t(u_t,w_t)\geq c \|u_t\|^\alpha_V +d \|u_t\|_V +d_0, \quad \forall \,  u_t \in \mathcal{K}_t, \; \forall \,  t \in T \,.
$$
\end{description}
\begin{theorem} \label{conv} Under the hypotheses $(H1)$-$(H6)$, there exists a solution to the problem  $VI(\psi_t,f_t,\mathcal{K}_t)$. Moreover, the family $\{u_t\}$ of solutions to the problem  $VI(\psi_t,f_t,\mathcal{K}_t)$ is bounded in $X$ and there exists a subnet of $\{u_t\}$ that converges weakly  in $X$ to  a solution of  the problem $VI(\psi,f,\mathcal{K})$. 
\end{theorem}
The hypotheses (H1) and (H2) are due to Glowinski \cite{Gl}  and describe the Mosco convergence \cite{At} of the family $\mathcal{K}_t$ to $\mathcal{K}$, whereas $f_t$ in (H4) is a standard approximation of the
linear functional $f$, for example, by numerical integration. The hypotheses (H5)-(H6) have been verified  in 
\cite{Ovcharova2015}.  

We apply Theorem \ref{conv} with $t=(\varepsilon, h, p)$ and show convergence for $\varepsilon \to 0^+$, and either $p\to \infty$ or $h\to 0$. We emphasize that in fact the weak convergence of $\{\mathbf{u}^\varepsilon_{hp}\}$ to $\mathbf{u}$ in $\mathbf{H}^{1/2}(\Gamma)$ can be replaced by the strong one. The complete proof of this result for the $h$-version of BEM can be found in \cite{Ovcharova2015}. 

\section{A-priori error estimate for the regularized problem} \label{sec:apriori}

 In this section we provide an a-priori error estimate for the $hp$-approximate solution $\mathbf{u}^{\varepsilon}_{hp}$ of the regularized boundary variational problem $(\mathcal{P}_\varepsilon)$ under the regularity assumptions of \cite{MS-1} for Signorini contact, i.e.~$\mathbf{u}_\varepsilon \in \mathbf{H}^{3/2}(\Gamma)$ and $g\in H^{3/2}(\Gamma_C)$. For our more general variational problem  $(\mathcal{P}_\varepsilon)$ we arrive at the same convergence order of $ \mathcal{O}(h^{1/4}p^{-1/4})$ as Maischak and Stephan in \cite{MS-1}.   
To this end,  we  present the following  C\'{e}a-Falk approximation lemma for the regularized problem $(\mathcal{P}_\varepsilon)$. 
\begin{lemma} Let $\mathbf{u}_\varepsilon \in \mathcal{K}^\Gamma$ be the solution of the problem $(\mathcal{P}_\varepsilon)$ and let  $\mathbf{u}^{\varepsilon}_{hp} \in \mathcal{K}_{hp}^\Gamma$ be the solution of the problem $(\mathcal{P}_{\varepsilon, hp})$. Assume that $\alpha_0<c_p$, $\mathbf{u}_\varepsilon$ $\in \mathbf{H}^{3/2}(\Gamma)$, $g\in H^{3/2}(\Gamma_C)$ and $P \mathbf{u}_\varepsilon -\mathbf{F} \in \mathbf{L}^2(\Gamma)$.  Then there exists a constant $c=c(\mathbf{u}_\varepsilon, g, \mathbf{f}, \mathbf{t})>0$, but independent of $h$ and $p$ such that
\begin{align*}
c&\|\mathbf{u}_\varepsilon-\mathbf{u}^\varepsilon_{ hp}\|^2_{\mathbf{H}^{1/2}(\Gamma)}  \leq  \|E_{hp}(\mathbf{u}_\varepsilon)\|^2_{\mathbf{H}^{-1/2}(\Gamma)} \\
& \quad +\inf _{v\in \mathcal{K}^\Gamma} 
\{ \|P \mathbf{u}_\varepsilon-\mathbf{F}\|_{\mathbf{L}^2(\Gamma)} \|\mathbf{u}^{\varepsilon}_{ hp}-\mathbf{v}\|_{\mathbf{L}^2(\Gamma)} + 
\langle D J_\varepsilon(\mathbf{u}_\varepsilon),\mathbf{v}-\mathbf{u}^\varepsilon_{hp} \rangle_{\Gamma_C} \} \\
&\quad +  \inf _{v_{hp}\in \mathcal{K}_{hp}^\Gamma} \left\{ \|\mathbf{u}_\varepsilon- \mathbf{v}_{ hp}\|_{\mathbf{H}^{1/2}(\Gamma)}^2 + \|P \mathbf{u}_\varepsilon-\mathbf{F}\|_{\mathbf{L}^2(\Gamma)} \|\mathbf{u}_\varepsilon-\mathbf{v}_{ hp}\|_{\mathbf{L}^2(\Gamma)} \right.\\
&\quad \left.+ 
\langle D J_\varepsilon(\mathbf{u}^{\varepsilon}_{ hp}),\mathbf{v}_{hp}-\mathbf{u}_{\varepsilon} \rangle_{\Gamma_C}
\right\}. 
\end{align*}
\end{lemma}
{\bf Proof} Using the definitions of $(\mathcal{P}_\varepsilon)$ and $(\mathcal{P}_{\varepsilon,hp})$, and estimates similar to  \cite[Theorem 3]{MS-1}, we obtain for all $v\in \mathcal{K}^\Gamma$,  $v_{hp}\in \mathcal{K}_{hp}^\Gamma$
\begin{align*}
 c_P &\|\mathbf{u}_\varepsilon-\mathbf{u}^\varepsilon_{ hp}\|^2_{\mathbf{H}^{1/2}(\Gamma)} \leq  \|E_{hp}(\mathbf{u}_\varepsilon)\|^2_{\mathbf{H}^{-1/2}(\Gamma)} + 
\|\mathbf{u}_\varepsilon- \mathbf{v}_{ hp}\|_{\mathbf{H}^{1/2}(\Gamma)}^2 \\
& + \|P \mathbf{u}_\varepsilon-\mathbf{F}\|_{\mathbf{L}^2(\Gamma)} \left( \|\mathbf{u}_\varepsilon-\mathbf{v}_{ hp}\|_{\mathbf{L}^2(\Gamma)}+   \|\mathbf{u}^{\varepsilon}_{ hp}-\mathbf{v}\|_{\mathbf{L}^2(\Gamma)} \right )  + {\rm D},
\end{align*}
where we abbreviate
\[
{\rm D}= \langle D J_\varepsilon(\mathbf{u}_\varepsilon),\mathbf{v}-\mathbf{u}_\varepsilon \rangle_{\Gamma_C} + \langle D J_\varepsilon(\mathbf{u}^{\varepsilon}_{ hp}),\mathbf{v}_{hp}-\mathbf{u}^\varepsilon_{hp} \rangle_{\Gamma_C}.
\]
To bound the term ${\rm D}$, we  use (\ref{ass1_th5}) and estimate as follows:
\begin{align*}
{\rm D} & =  \langle D J_\varepsilon(\mathbf{u}_\varepsilon),\mathbf{v}-\mathbf{u}^\varepsilon_{hp} \rangle_{\Gamma_C} + \langle D J_\varepsilon(\mathbf{u}^{\varepsilon}_{ hp}), \mathbf{v}_{hp}-\mathbf{u}_{\varepsilon} \rangle_{\Gamma_C} \\
& \quad + 
\langle D J_\varepsilon(\mathbf{u}_\varepsilon)-  D J_\varepsilon(\mathbf{u}^{\varepsilon}_{ hp}) ,\mathbf{u}^\varepsilon_{hp} - \mathbf{u}_\varepsilon \rangle_{\Gamma_C} \\
 & \leq  \langle D J_\varepsilon(\mathbf{u}_\varepsilon),\mathbf{v}-\mathbf{u}^\varepsilon_{hp} \rangle_{\Gamma_C} + \langle D J_\varepsilon(\mathbf{u}^{\varepsilon}_{ hp}), \mathbf{v}_{hp}-\mathbf{u}_{\varepsilon} \rangle_{\Gamma_C} + \alpha_0 \| \mathbf{u}_{\varepsilon} - \mathbf{u}^\varepsilon_{ hp}\|^2_{\mathcal{V}}. 
\end{align*}
Therefore, since $\alpha_0 <c_P$, see Theorem \ref{theouniq}, we obtain the assertion. \qed
\begin{theorem}
Let $\mathbf{u}_\varepsilon \in \mathcal{K}^\Gamma$ be the solution of the problem $(\mathcal{P}_\varepsilon)$ and let $\mathbf{u}^{\varepsilon}_ {hp} \in \mathcal{K}_{hp}^\Gamma$ be the solution of the problem $(\mathcal{P}_{\varepsilon, hp})$. Assume that $\alpha_0<c_p$, $\mathbf{u}_\varepsilon$ $\in \mathbf{H}^{3/2}(\Gamma)$, $g\in H^{3/2}(\Gamma_C)$ and $P \mathbf{u}_\varepsilon -\mathbf{F} \in \mathbf{L}^2(\Gamma)$.  Then there exists a constant  $c=c(\mathbf{u}_\varepsilon, g, \mathbf{f}, \mathbf{t})>0$, but independent of $h$ and $p$ such that
\begin{equation} \label{a-priori_est}
\|\mathbf{u}_\varepsilon-\mathbf{u}^{\varepsilon}_{hp}\|_{\mathbf{H}^{1/2}(\Gamma)} \leq c h^{1/4}p^{-1/4}.
\end{equation}
\end{theorem} 
{\bf Proof} 
Taking into account the estimates obtained by Maischak and Stephan in \cite[Theorem 3]{MS-1} for  the  consistency error, the approximation error, and for $\|E_{hp}{\mathbf{u}}\|_{\mathbf{H}^{-1/2}(\Gamma)}$, 
we only need to estimate 
\begin{equation}\label{bound1}
\inf\{\langle D J_\varepsilon(\mathbf{u}_\varepsilon),\mathbf{v}-\mathbf{u}^\varepsilon_{hp} \rangle_{\Gamma_C} \, : \,  \mathbf{v} \in \mathcal{K}^\Gamma \}
\end{equation}
and
\begin{equation}\label{bound2}
\inf \{\langle D J_\varepsilon(\mathbf{u}^{\varepsilon}_{ hp}),\mathbf{v}_{hp}-\mathbf{u}_{\varepsilon} \rangle_{\Gamma_C} \, : \,  \mathbf{v}_{hp}\in \mathcal{K}_{hp}^\Gamma\}.
\end{equation}
To estimate (\ref{bound2}) we take $\mathbf{v}_{hp}:= I_{hp} \mathbf{u}_\varepsilon \in \mathcal{K}_{hp}^{\Gamma}$, the interpolate of $\mathbf{u}_\varepsilon\in \mathbf{H}^{3/2}(\Gamma)\subset C^0(\Gamma)$.
From (\ref{ass_3}),
\begin{align}
 \langle D J_\varepsilon(\mathbf{u}^{\varepsilon}_{ hp}),\mathbf{v}_{hp}-\mathbf{u}_{\varepsilon} \rangle_{\Gamma_C} 
& =     \int_{\Gamma_C} S_x(u^\varepsilon_{hp, n}, \varepsilon)(v_{hp, n}- u_{\varepsilon, n}) \, ds \nonumber\\
  & \leq   c \left(1+\|u^\varepsilon_{ hp, n}\|_{L^2(\Gamma_C)}\right) \| u_{\varepsilon, n}- v_{hp, n}\|_{L^2(\Gamma_C)}.  \label{bound3}
\end{align}
By \cite[Theorems 4.2 and 4.5] {Bernardi}  and by the real interpolation between $H^1(\Gamma)$ and $L^2(\Gamma)$ there exists a constant $C_1>0$ such that
\begin{equation} \label{bound4}
\|u_{\varepsilon, n} - v_{hp, n}\|_{H^{1/2}(\Gamma)} \leq C_1 h^1 p^{-1} \, \|\mathbf{u}_{\varepsilon}\|_{H^{3/2}(\Gamma)}.
\end{equation}
To estimate (\ref{bound1}) similarly as \cite{MS-1} we define 
$\mathbf{v}^* \in \mathcal{K}^\Gamma \cap \mathbf{H}^1(\Gamma)$ by 
\[
\mathbf{v}^* :=\left\{ \begin{array}{lc}
\mathbf{u}^\varepsilon_{hp,t}+[g+ \inf \{u^\varepsilon_{hp,n}-g_{hp},0\}] \mathbf{n} \;&  \mbox{on}\; \Gamma_C \\
0 \;&  \mbox{on}\; \Gamma_D \\
\gamma_N \mathbf{u}^\varepsilon_{hp} \; &  \mbox{on}\; \Gamma_N,
\end{array} 
\right.
\]
where $g_{hp}:=I_{hp} g$ is the interpolate of  the gap function $g$, and $\gamma _N$ is the trace map onto $\Gamma_N$.  
\\
Analogously to (\ref{bound3}), we have
\begin{align}
\langle D J_\varepsilon(\mathbf{u}_\varepsilon),\mathbf{v}^*-\mathbf{u}^\varepsilon_{hp} \rangle_{\Gamma_C} & =   \int_{\Gamma_C} S_x(u_{\varepsilon, n}, \varepsilon)(v^*_n- u^\varepsilon_{hp, n}) \, ds\nonumber\\
  & \leq  c \left(1+\|u_{\varepsilon, n}\|_{L^2(\Gamma_C)}\right) \|v^*_n- u^\varepsilon_{hp, n}\|_{L^2(\Gamma_C)}. \label{bound5}
\end{align}
The elaborate analysis in \cite{MS-1}, see the proof of Theorem 3, estimates (31) - (33),  gives
\begin{equation}\label{bound6}
\|v^*_n-u^\varepsilon_{hp,n}\|_{L^2(\Gamma_C)} \leq C_2 h^{1/2}p^{-1/2}\left( \|g\|_{H^{1/2}(\Gamma_C)}+ \|\mathbf{u}^{\varepsilon}_{hp}\|_{H^{1/2}(\Gamma_C)} \right).
\end{equation}
\\
Finally, combining the error estimates for the interpolation  (\ref{bound4}) and the consistency (\ref{bound6}) with (\ref{bound3}) and  (\ref{bound5}), respectively, and taking into account the boundedness of $\|\mathbf{u}^{\varepsilon}_{hp}\|$ in $H^{1/2}(\Gamma_C)$ (see Theorem \ref{conv}),  we prove  the wished bound for (\ref{bound1}) and (\ref{bound2}). \qed
\begin{rem} Let $\mathbf{u} \in \mathcal{K}^\Gamma$ be the solution of the problem $(\mathcal{P})$ and $\mathbf{u}_\varepsilon \in \mathcal{K}^\Gamma$ be the solution of the problem $(\mathcal{P}_\varepsilon)$.
Taking $\mathbf{v}=\mathbf{u}_\varepsilon$ in   (\ref{BdHemVar}) and $\mathbf{v}=\mathbf{u}$ in  (\ref{reg}), and adding the two resulting inequalities yields
\[
\langle P \mathbf{u}-P \mathbf{u}_\varepsilon, \mathbf{u}_\varepsilon- \mathbf{u} \rangle _{\Gamma_0}+ \varphi(\mathbf{u}, \mathbf{u}_\varepsilon) + \varphi_\varepsilon( \mathbf{u}_\varepsilon, \mathbf{u}) \geq 0.
\]
Hence,
\[
c_P\|\mathbf{u}_\varepsilon-\mathbf{u}\|^2_{\mathcal {V}}\leq \langle P \mathbf{u}_\varepsilon-P \mathbf{u}, \mathbf{u}_\varepsilon- \mathbf{u} \rangle _{\Gamma_0} \leq \varphi(\mathbf{u}, \mathbf{u}_\varepsilon) + \varphi_\varepsilon( \mathbf{u}_\varepsilon, \mathbf{u}). 
\]
We emphasize that, in general, the functional $\varphi_\varepsilon$ does not approximate the functional $\varphi$ arbitrarily close as $\varepsilon \to 0^+$. Nevertheless, we have convergence of the sequence $\{\mathbf{u}_\varepsilon\}$ of the solutions of the regularized problem ($\mathcal{P}_\varepsilon$) to the solution $\mathbf{u}$ of the boundary hemivariational problem ($\mathcal{P}$). In particular, for any $\mathbf{v}_\varepsilon \rightharpoonup \mathbf{v}$  it holds
\[
\limsup_{\varepsilon\to 0} \, \varphi_\varepsilon (\mathbf{v}_\varepsilon, \mathbf{v}) \leq \varphi(\mathbf{v}, \mathbf{v}) \quad \forall \mathbf{v} \in \mathcal{K}^\Gamma,
\]
which according to Theorem \ref{conv} is a sufficient condition for convergence. \\
However, for convergence rates in $\varepsilon$,
\[
\varphi(\mathbf{u}, \mathbf{u}_\varepsilon) + \varphi_\varepsilon( \mathbf{u}_\varepsilon, \mathbf{u}) \leq c \, \varepsilon^\beta \|\mathbf{u}_\varepsilon-\mathbf{u}\|_{\mathcal{V}}
\]
for some positive constants $c$ and $\beta$ is needed. An estimate of this form can be found in \cite{Ovcharova2012} for the approximation of the absolute value function.
\end{rem}

\section{A-posteriori error estimate for the regularized problem}
\label{sec:aposteriori}
To be able to split the approximation error into the discretization error of a simpler variational equation and contributions arising from the constraints on $\Gamma_C$ we introduce the regularized mixed formulation \eqref{eq:MixedVarEq}-\eqref{eq:mixedprobConstraint}, which is equivalent to the regularized problem $(\mathcal{P}_\varepsilon)$. \\
Find $(\mathbf{u}^\varepsilon,\lambda^\varepsilon) \in \mathcal{V} \times M(\mathbf{u}^\varepsilon)$ such that
\begin{subequations} \label{eq:mixedprob}
\begin{alignat}{2} 
 \left\langle P \mathbf{u}^\varepsilon, \mathbf{v} \right\rangle_{\Gamma_0}+ \left\langle \lambda^\varepsilon, v_n \right\rangle_{\Gamma_C} & = \langle \mathbf{F}, \mathbf{v} \rangle_{\Gamma_0} & \quad & \forall \,  \mathbf{v} \in \mathcal{V} \label{eq:MixedVarEq}\\
 \left\langle \mu -\lambda^\varepsilon, u^\varepsilon_n - g \right\rangle_{\Gamma_C} & \leq 0 & \quad & \forall \,  \mu \in M(\mathbf{u}^\varepsilon) \label{eq:mixedprobConstraint}
\end{alignat}
\end{subequations}
with the set of admissible Lagrange multipliers
\begin{align*}
 M(\mathbf{u}^\varepsilon):=\left\{ \mu \in X^*:\ \left\langle \mu,\eta \right\rangle_{\Gamma_C} \geq \left\langle DJ_\varepsilon(\mathbf{u}^\varepsilon),\eta \right\rangle_{\Gamma_C} \ \forall \,  \eta \in X,\ \eta \geq 0 \text{ a.e.~on } \Gamma_C \right\} 
\end{align*}
where $X=\{w \mid \exists\, \mathbf{v} \in \mathcal{V},\ v_n|_{\Gamma_C} = w\}\subset H^{1/2}(\Gamma_C)$ and $X^*$ its dual space.
\begin{lemma}
 \begin{enumerate}
  \item Let $\mathbf{u}^\varepsilon$ solve the regularized problem  $(\mathcal{P}_\varepsilon)$, then there exists a $\lambda^\varepsilon \in  M(\mathbf{u}^\varepsilon)$ such that $(\mathbf{u}^\varepsilon,\lambda^\varepsilon)$ solves \eqref{eq:mixedprob}. \\[0.01cm]
  \item Let $(\mathbf{u}^\varepsilon,\lambda^\varepsilon)$ solve \eqref{eq:mixedprob}, then $\mathbf{u}^\varepsilon$ solves $(\mathcal{P}_\varepsilon)$.
 \end{enumerate}
\end{lemma}
{\bf Proof}\, 
 \begin{enumerate}
  \item Define $\lambda^\varepsilon \in X^*$ by
\begin{align} \label{eq1_mixed}
 \left\langle \lambda^\varepsilon, v_n \right\rangle_{\Gamma_C}  = \langle \mathbf{F}, \mathbf{v} \rangle_{\Gamma_0} - \left\langle P \mathbf{u}^\varepsilon, \mathbf{v} \right\rangle_{\Gamma_0} \quad \forall \,  v \in \mathcal{V}.
\end{align}
Then, by \eqref{reg}, we obtain for all $\mathbf{v} \in \mathcal{K}^\Gamma$ with $v_n|_{\Gamma_C}=u^\varepsilon_n|_{\Gamma_C} -\eta$ for $0 \leq \eta \in  X$
\begin{gather*}
 \left\langle \lambda^\varepsilon, v_n -u^\varepsilon_n \right\rangle_{\Gamma_C}  = \langle \mathbf{F}, \mathbf{v}-\mathbf{u}^\varepsilon \rangle_{\Gamma_0} - \left\langle P \mathbf{u}^\varepsilon, \mathbf{v} -\mathbf{u}^\varepsilon \right\rangle_{\Gamma_0} \leq \left\langle DJ_\varepsilon(\mathbf{u}^\varepsilon),v_n -u^\varepsilon_n  \right\rangle_{\Gamma_C} \\  \Rightarrow \quad 
 \left\langle \lambda^\varepsilon, -\eta \right\rangle_{\Gamma_C} \leq \left\langle DJ_\varepsilon(\mathbf{u}^\varepsilon),-\eta  \right\rangle_{\Gamma_C},
\end{gather*}
i.e.~$\lambda^\varepsilon \in  M(\mathbf{u}^\varepsilon)$. Analogously, we obtain for $v_n|_{\Gamma_C}=g$, $v_n|_{\Gamma_C}=2u^\varepsilon_n|_{\Gamma_C} -g$
\begin{gather*}
\left\langle \lambda^\varepsilon, g -u^\varepsilon_n \right\rangle_{\Gamma_C} \leq \left\langle DJ_\varepsilon(\mathbf{u}^\varepsilon),g -u^\varepsilon_n  \right\rangle_{\Gamma_C} \  \text{and} \ 
\left\langle \lambda^\varepsilon, u^\varepsilon_n -g \right\rangle_{\Gamma_C} \leq \left\langle DJ_\varepsilon(\mathbf{u}^\varepsilon),u^\varepsilon_n -g  \right\rangle_{\Gamma_C},
\end{gather*}
which implies 
\[
 \left\langle \lambda^\varepsilon, u^\varepsilon_n -g \right\rangle_{\Gamma_C} = \left\langle DJ_\varepsilon(\mathbf{u}^\varepsilon),u^\varepsilon_n -g \right\rangle_{\Gamma_C}.
\]
It remains to show \eqref{eq:mixedprobConstraint}. From the previous equation we obtain for $\mu \in M(\mathbf{u}^\varepsilon)$
\begin{align*}
 \left\langle \mu -\lambda^\varepsilon, u^\varepsilon_n - g \right\rangle_{\Gamma_C} &= 
 \left\langle \mu , u^\varepsilon_n - g \right\rangle_{\Gamma_C} - \left\langle DJ_\varepsilon(\mathbf{u}^\varepsilon),u^\varepsilon_n -g \right\rangle_{\Gamma_C} \\
 &\leq \left\langle DJ_\varepsilon(\mathbf{u}^\varepsilon),u^\varepsilon_n -g \right\rangle_{\Gamma_C} - \left\langle DJ_\varepsilon(\mathbf{u}^\varepsilon),u^\varepsilon_n -g \right\rangle_{\Gamma_C} =0.
\end{align*}

\item There holds trivially $ DJ_\varepsilon(\mathbf{u}^\varepsilon)$, $2\lambda -  DJ_\varepsilon(\mathbf{u}^\varepsilon) \in M(\mathbf{u}^\varepsilon)$, and thus \eqref{eq:mixedprobConstraint} yields
\begin{gather*}
\left\langle -\lambda^\varepsilon, u^\varepsilon_n -g \right\rangle_{\Gamma_C} \leq \left\langle -DJ_\varepsilon(\mathbf{u}^\varepsilon),u^\varepsilon_n -g  \right\rangle_{\Gamma_C},
 \ \left\langle \lambda^\varepsilon, g -u^\varepsilon_n \right\rangle_{\Gamma_C} \leq \left\langle DJ_\varepsilon(\mathbf{u}^\varepsilon),g -u^\varepsilon_n  \right\rangle_{\Gamma_C} 
 \\  \Rightarrow \quad
 \left\langle \lambda^\varepsilon, u^\varepsilon_n -g \right\rangle_{\Gamma_C} = \left\langle DJ_\varepsilon(\mathbf{u}^\varepsilon),u^\varepsilon_n -g \right\rangle_{\Gamma_C}.
\end{gather*}
For $\mathbf{v} \in \mathcal{K}^\Gamma$, i.e.~$v_n|_{\Gamma_C}=g -\eta$ with $0 \leq \eta \in X$, we therefore obtain
\begin{align*}
 \left\langle \lambda^\varepsilon, v_n -u^\varepsilon_n \right\rangle_{\Gamma_C}  &= \left\langle \lambda^\varepsilon, g -u^\varepsilon_n \right\rangle_{\Gamma_C}  +
 \left\langle \lambda^\varepsilon,  -\eta  \right\rangle_{\Gamma_C}  = \left\langle DJ_\varepsilon(\mathbf{u}^\varepsilon), g -u^\varepsilon_n \right\rangle_{\Gamma_C}  + \left\langle \lambda^\varepsilon,  -\eta  \right\rangle_{\Gamma_C} \\
 &\leq \left\langle DJ_\varepsilon(\mathbf{u}^\varepsilon), g -u^\varepsilon_n \right\rangle_{\Gamma_C} + \left\langle DJ_\varepsilon(\mathbf{u}^\varepsilon), -\eta \right\rangle_{\Gamma_C} = \left\langle DJ_\varepsilon(\mathbf{u}^\varepsilon), v_n -u^\varepsilon_n \right\rangle_{\Gamma_C}.
\end{align*}
It remains to show that $u^\varepsilon_n -g \leq 0$. Assume there exists a measurable set in which $u^\varepsilon_n - g>0$. Then, from \eqref{eq:mixedprobConstraint} we obtain for $\mu=DJ_\varepsilon(\mathbf{u}^\varepsilon) + \chi_{u^\varepsilon_n - g >0} \in M(\mathbf{u}^\varepsilon)$ with indicator function $\chi$
\begin{align*}
 0 \geq \left\langle \mu -\lambda^\varepsilon, u^\varepsilon_n - g \right\rangle_{\Gamma_C}
 =\left\langle \mu -DJ_\varepsilon(\mathbf{u}^\varepsilon), u^\varepsilon_n - g \right\rangle_{\Gamma_C}
 =\left\langle \chi_{u^\varepsilon_n - g >0}, u^\varepsilon_n - g \right\rangle_{\Gamma_C} >0,
\end{align*}
i.e.~$\mathbf{u}^\varepsilon \in \mathcal{K}^\Gamma$.
\qed
\end{enumerate}
Given the discrete solution $\mathbf{u}_{hp}^\varepsilon \in \mathcal{K}_{hp}^\Gamma$ to $(\mathcal{P}_{\varepsilon,hp})$, we reconstruct $\lambda_{hp}^\varepsilon \in \operatorname{span}\left\{ \psi_i \right\}_{i=1}^M $ such that
\begin{align} \label{eq:reconstruct_lambda_h}
  \left\langle \lambda_{hp}^\varepsilon, v_n \right\rangle_{\Gamma_C}  = \langle \mathbf{F}, \mathbf{v} \rangle_{\Gamma_0} - \left\langle P_{hp} \mathbf{u}_{hp}^\varepsilon, \mathbf{v} \right\rangle_{\Gamma_0} \quad \forall \,  v \in \mathcal{V}_{hp}
\end{align}
by solving a potentially over-constrained system of linear equations for an arbitrary choice of basis $\left\{ \psi \right\}$.
Following the Braess trick \cite{braess2005posteriori} as e.g.~in \cite{Banz2015}, we define the auxiliary problem
\begin{align} \label{eq:auxProb}
  \mathbf{z} \in \mathcal{V} : \quad  \left\langle P \mathbf{z}, \mathbf{v} \right\rangle_{\Gamma_0} = \langle \mathbf{F}, \mathbf{v} \rangle_{\Gamma_0} - \left\langle \lambda_{hp}^\varepsilon, v_n \right\rangle_{\Gamma_C} \quad  \forall \,  \mathbf{v} \in \mathcal{V}.
\end{align}
Subtracting \eqref{eq:MixedVarEq} and \eqref{eq:auxProb} yields
\begin{align} \label{eq:exchange_u_to_z} 
\left\langle P(\mathbf{u}^\varepsilon -\mathbf{z}), \mathbf{v} \right\rangle_{\Gamma_0} =  \left\langle \lambda_{hp}^\varepsilon - \lambda^\varepsilon, v_n \right\rangle_{\Gamma_C} \quad  \forall \, \mathbf{v} \in \mathcal{V} 
\end{align}
and additionally with the continuous inf-sup condition \cite[Theorem 3.2.1]{Chernov} this yields (see \cite{Banz2015})
\begin{align} \label{eq:apostLambda}
 \left\| \lambda_{hp}^\varepsilon - \lambda^\varepsilon \right\|_{X^*} \leq \frac{C}{\beta} \left\| \mathbf{u}^\varepsilon - \mathbf{z} \right\|_{\mathcal{V}} \leq \frac{C}{\beta} \left\| \mathbf{u}^\varepsilon - \mathbf{u}_{hp}^\varepsilon \right\|_{\mathcal{V}} + \frac{C}{\beta} \left\| \mathbf{u}_{hp}^\varepsilon - \mathbf{z} \right\|_{\mathcal{V}}
 \end{align}
with inf-sup constant $\beta>0$. See \cite[Theorem 3.2.1]{Chernov} for a proof of the inf-sup condition for the difficult case when $\bar{\Gamma}_C \cap \bar{\Gamma}_D = \emptyset$, i.e.~$X^*=\tilde{H}^{-1/2}(\Gamma_C)$.

\begin{theorem} \label{thm:aposteriori}
Under the assumption \eqref{uniqreg} and if $S_x(\cdot, \varepsilon)$ is Lipschitz continuous, then there exists a constant $C$ independent of $h$ and $p$ such that for $\varsigma>0$ arbitrary
\begin{align*}
 \left(c_P - \alpha_0 -4\varsigma \right)\|\mathbf{u}^\varepsilon - \mathbf{u}_{hp}^\varepsilon \|_{\mathcal{V}}^2 \leq &  \left(\frac{C}{\varsigma} + 1 \right) \|\mathbf{z} - \mathbf{u}_{hp}^\varepsilon \|^2_{\mathcal{V}}  + \frac{1}{4\varsigma}\left\| ( \lambda_{hp}^\varepsilon -DJ_{\varepsilon}(\mathbf{u}_{hp}^\varepsilon) )^- \right\|_{ X^*}^2 
 \\
 &+ C\left(\frac{1}{\varsigma} + \frac{1}{\beta^2}+\frac{1}{\varsigma \beta^2}\right) \left\|(u_{hp,n}^\varepsilon -g)^+ \right\|^2_{ X} \\
 & -  \left\langle   (\lambda_{hp}^\varepsilon -  DJ_{\varepsilon}(\mathbf{u}_{hp}^\varepsilon))^+ , (u_{hp,n}^\varepsilon -g)^- \right\rangle_{\Gamma_C}
\end{align*}
with $\lambda_{hp}^\varepsilon$ satisfying \eqref{eq:reconstruct_lambda_h}.
\end{theorem}
{\bf Proof}\, 
From the Lipschitz continuity it follows that 
\begin{align*}
 \left\langle  DJ_\varepsilon(\mathbf{u}_{hp}^\varepsilon) - DJ_\varepsilon(\mathbf{u}^\varepsilon), v_n  \right\rangle_{\Gamma_C} &= \int_{\Gamma_C}  \left( S_x({u}_{hp,n}^\varepsilon) - S_x({u}_{n}^\varepsilon) \right) v_n \, ds\\
 &\leq Lip \|{u}_{hp,n}^\varepsilon - {u}_{n}^\varepsilon \|_{L^2(\Gamma_C)} \|{v}_{n} \|_{L^2(\Gamma_C)}
\end{align*}
and thus, by Young's inequality ($\varsigma>0$ arbitrary), 
\begin{align} 
\left\langle DJ_{\varepsilon}(\mathbf{u}_{hp}^\varepsilon) -  DJ_\varepsilon(\mathbf{u}^\varepsilon)  , (u_{hp,n}^\varepsilon -g)^+ \right\rangle_{\Gamma_C} & \leq  \varsigma \|\mathbf{u}_{hp}^\varepsilon - \mathbf{u}^\varepsilon \|_{L^2(\Gamma_C)}^2  \nonumber \\
& \quad + \frac{C^2}{4\varsigma} \|(u_{hp,n}^\varepsilon -g)^+ \|_{L^2(\Gamma_C)}^2. \label{eq:apostHelp}
\end{align}
From the conformity $\mathcal{V}_{hp} \subset \mathcal{V}$ we obtain with \eqref{eq:exchange_u_to_z}
\begin{align*}
 c_P \|\mathbf{u}^\varepsilon - \mathbf{u}_{hp}^\varepsilon \|_{\mathcal{V}}^2 &\leq \left\langle P(\mathbf{u}^\varepsilon -\mathbf{u}_{hp}^\varepsilon), \mathbf{u}^\varepsilon -\mathbf{u}_{hp}^\varepsilon \right\rangle_{\Gamma_0}\\
 &= \left\langle P(\mathbf{z}-\mathbf{u}_{hp}^\varepsilon), \mathbf{u}^\varepsilon -\mathbf{u}_{hp}^\varepsilon \right\rangle_{\Gamma_0} + \left\langle \lambda_{hp}^\varepsilon - \lambda^\varepsilon, u^\varepsilon_n -u_{hp,n}^\varepsilon \right\rangle_{\Gamma_C} \\
 &\leq C_P \|\mathbf{z} - \mathbf{u}_{hp}^\varepsilon \|_{\mathcal{V}} \|\mathbf{u}^\varepsilon - \mathbf{u}_{hp}^\varepsilon \|_{\mathcal{V}} + \left\langle \lambda_{hp}^\varepsilon - \lambda^\varepsilon, u^\varepsilon_n -u_{hp,n}^\varepsilon \right\rangle_{\Gamma_C}.
\end{align*}
From complementarity $ \left\langle \lambda^\varepsilon, u^\varepsilon_n -g \right\rangle_{\Gamma_C} = \left\langle DJ_\varepsilon(\mathbf{u}^\varepsilon),u^\varepsilon_n -g \right\rangle_{\Gamma_C}$, $\lambda^\varepsilon \in M(\mathbf{u}^\varepsilon)$, $v=v^+ + v^- = \max\{0,v\} + \min\{0,v\}$, $u^\varepsilon_n - g \leq 0$ a.e.~on $\Gamma_C$ we obtain
\begin{align*}
  & \left\langle \lambda_{hp}^\varepsilon - \lambda^\varepsilon, u^\varepsilon_n -u_{hp,n}^\varepsilon \right\rangle_{\Gamma_C} = 
   \left\langle \lambda_{hp}^\varepsilon - DJ_{\varepsilon}(\mathbf{u}^\varepsilon_{hp}), u^\varepsilon_n -u_{hp,n}^\varepsilon \right\rangle_{\Gamma_C}
\\  &\quad  +    \left\langle  DJ_{\varepsilon}(\mathbf{u}^\varepsilon_{hp})- \lambda^\varepsilon, u^\varepsilon_n -g \right\rangle_{\Gamma_C}
+    \left\langle  DJ_{\varepsilon}(\mathbf{u}^\varepsilon_{hp})- \lambda^\varepsilon, g -u_{hp,n}^\varepsilon \right\rangle_{\Gamma_C} \\
 & = 
   \left\langle \lambda_{hp}^\varepsilon - DJ_{\varepsilon}(\mathbf{u}^\varepsilon_{hp}), u^\varepsilon_n -u_{hp,n}^\varepsilon \right\rangle_{\Gamma_C}
+    \left\langle  DJ_{\varepsilon}(\mathbf{u}^\varepsilon_{hp})- DJ_\varepsilon(\mathbf{u}^\varepsilon), u^\varepsilon_n -g \right\rangle_{\Gamma_C}
\\  &\quad+    \left\langle  DJ_{\varepsilon}(\mathbf{u}^\varepsilon_{hp})- \lambda^\varepsilon, g -u_{hp,n}^\varepsilon \right\rangle_{\Gamma_C}\\
 & = 
   \left\langle \lambda_{hp}^\varepsilon - DJ_{\varepsilon}(\mathbf{u}^\varepsilon_{hp}), u^\varepsilon_n -u_{hp,n}^\varepsilon \right\rangle_{\Gamma_C}
+    \left\langle  DJ_{\varepsilon}(\mathbf{u}^\varepsilon_{hp})- DJ_\varepsilon(\mathbf{u}^\varepsilon), u^\varepsilon_n -u_{hp,n}^\varepsilon  \right\rangle_{\Gamma_C}
\\  &\quad+    \left\langle \lambda^\varepsilon -  DJ_\varepsilon(\mathbf{u}^\varepsilon) , u_{hp,n}^\varepsilon -g \right\rangle_{\Gamma_C}\\
& \leq  \left\langle (\lambda_{hp}^\varepsilon - DJ_{\varepsilon}(\mathbf{u}^\varepsilon_{hp}))^+, u^\varepsilon_n -g+g -u_{hp,n}^\varepsilon \right\rangle_{\Gamma_C} + \left\langle (\lambda_{hp}^\varepsilon - DJ_{\varepsilon}(\mathbf{u}^\varepsilon_{hp}))^-, u^\varepsilon_n -u_{hp,n}^\varepsilon \right\rangle_{\Gamma_C}\\
&\quad +    \left\langle  DJ_{\varepsilon}(\mathbf{u}^\varepsilon_{hp})- DJ_\varepsilon(\mathbf{u}^\varepsilon), u^\varepsilon_n -u_{hp,n}^\varepsilon  \right\rangle_{\Gamma_C}
+    \left\langle \lambda^\varepsilon -  DJ_\varepsilon(\mathbf{u}^\varepsilon) , (u_{hp,n}^\varepsilon -g)^+ \right\rangle_{\Gamma_C}
\\
& \leq  \left\langle (\lambda_{hp}^\varepsilon - DJ_{\varepsilon}(\mathbf{u}^\varepsilon_{hp}))^+, g -u_{hp,n}^\varepsilon \right\rangle_{\Gamma_C} +    \left\langle \lambda^\varepsilon -  DJ_\varepsilon(\mathbf{u}^\varepsilon) , (u_{hp,n}^\varepsilon -g)^+ \right\rangle_{\Gamma_C}\\
&\quad +    \left\langle  DJ_{\varepsilon}(\mathbf{u}^\varepsilon_{hp})- DJ_\varepsilon(\mathbf{u}^\varepsilon), u^\varepsilon_n -u_{hp,n}^\varepsilon  \right\rangle_{\Gamma_C} + \left\langle (\lambda_{hp}^\varepsilon - DJ_{\varepsilon}(\mathbf{u}^\varepsilon_{hp}))^-, u^\varepsilon_n -u_{hp,n}^\varepsilon \right\rangle_{\Gamma_C}.
\end{align*}
Since
\begin{align*}
&\left\langle (\lambda_{hp}^\varepsilon - DJ_{\varepsilon}(\mathbf{u}^\varepsilon_{hp}))^+, g -u_{hp,n}^\varepsilon \right\rangle_{\Gamma_C}+ \left\langle \lambda^\varepsilon -  DJ_\varepsilon(\mathbf{u}^\varepsilon) , (u_{hp,n}^\varepsilon -g)^+ \right\rangle_{\Gamma_C} \\
&=
- \left\langle (\lambda_{hp}^\varepsilon - DJ_{\varepsilon}(\mathbf{u}^\varepsilon_{hp}))^+, (u_{hp,n}^\varepsilon -g)^+ +(u_{hp,n}^\varepsilon -g)^-  \right\rangle_{\Gamma_C} 
\\& \quad +  \left\langle \lambda^\varepsilon -  DJ_\varepsilon(\mathbf{u}^\varepsilon)  - (\lambda_{hp}^\varepsilon -  DJ_{\varepsilon}(\mathbf{u}_{hp}^\varepsilon)) +  (\lambda_{hp}^\varepsilon -  DJ_{\varepsilon}(\mathbf{u}_{hp}^\varepsilon))^+ , (u_{hp,n}^\varepsilon -g)^+ \right\rangle_{\Gamma_C} \\
& \quad + \left\langle (\lambda_{hp}^\varepsilon -  DJ_{\varepsilon}(\mathbf{u}_{hp}^\varepsilon))^-, (u_{hp,n}^\varepsilon -g)^+ \right\rangle_{\Gamma_C} \\
&\leq \left\langle \lambda^\varepsilon   - \lambda_{hp}^\varepsilon , (u_{hp,n}^\varepsilon -g)^+ \right\rangle_{\Gamma_C}  +  \left\langle DJ_{\varepsilon}(\mathbf{u}_{hp}^\varepsilon) -  DJ_\varepsilon(\mathbf{u}^\varepsilon)  , (u_{hp,n}^\varepsilon -g)^+ \right\rangle_{\Gamma_C} \\
& \quad -  \left\langle   (\lambda_{hp}^\varepsilon -  DJ_{\varepsilon}(\mathbf{u}_{hp}^\varepsilon))^+ , (u_{hp,n}^\varepsilon -g)^- \right\rangle_{\Gamma_C},
\end{align*}
application of the Cauchy-Schwarz inequality, Young's inequality with $\varsigma >0$, \eqref{eq:apostLambda}, \eqref{eq:apostHelp} and \eqref{ass1_th5} yields
\begin{align*}
 \left(c_P - \alpha_0 -4\varsigma \right)\|\mathbf{u}^\varepsilon - \mathbf{u}_{hp}^\varepsilon \|_{\mathcal{V}}^2 \leq &  \left(\frac{C}{\varsigma} + 1 \right) \|\mathbf{z} - \mathbf{u}_{hp}^\varepsilon \|^2_{\mathcal{V}}  + \frac{1}{4\varsigma}\left\| ( \lambda_{hp}^\varepsilon -DJ_{\varepsilon}(\mathbf{u}_{hp}^\varepsilon) )^- \right\|_{ X^*}^2 
 \\
 &+ C\left(\frac{1}{\varsigma} + \frac{1}{\beta^2}+\frac{1}{\varsigma \beta^2}\right) \left\|(u_{hp,n}^\varepsilon -g)^+ \right\|^2_{ X} \\
 & -  \left\langle   (\lambda_{hp}^\varepsilon -  DJ_{\varepsilon}(\mathbf{u}_{hp}^\varepsilon))^+ , (u_{hp,n}^\varepsilon -g)^- \right\rangle_{\Gamma_C}.  \qquad \qed
\end{align*}

The a-posteriori error estimate decomposes into the discretization error of a variational equality $\|\mathbf{z} - \mathbf{u}_{hp}^\varepsilon \|^2_{\mathcal{V}}$, which can be further estimated by e.g.~residual error estimates \cite{cc96} or bubble error estimates, e.g.~\cite{Banz2015}, and violation of the consistency condition $\left\| ( \lambda_{hp}^\varepsilon -DJ_{\varepsilon}(\mathbf{u}_{hp}^\varepsilon) )^- \right\|_{ X^*}^2$, violation of the non-penetration condition $ \left\|(u_{hp,n}^\varepsilon -g)^+ \right\|^2_{ X}$ and violation of the complementarity condition $-  \left\langle   (\lambda_{hp}^\varepsilon -  DJ_{\varepsilon}(\mathbf{u}_{hp}^\varepsilon))^+ , (u_{hp,n}^\varepsilon -g)^- \right\rangle_{\Gamma_C}$. Localizing an approximation of the global a-posteriori error estimate gives rise to the following solve-mark-refine algorithm for $hp$-adaptivity.

 \begin{alg} \label{alg:adaptive}
\vspace*{2mm} \hrule
\vspace*{2mm} (Solve-mark-refine algorithm for $hp$-adaptivity) \vspace*{2mm}
\hrule
\begin{enumerate}
	\item Choose initial discretization $\mathcal{T}_{h}$ and $p$, steering parameters $\theta \in (0,1)$ and $\delta \in (0,1)$.
  \item  For $k=0,1,2,\ldots$ do	
 \begin{enumerate} 
  \item solve regularized discrete problem ($\mathcal{P}_{varepsilon, hp}$).
  \item compute discrete Lagrange multiplier by \eqref{eq:reconstruct_lambda_h}
  \item compute local error indicators $\Xi^2$ to current solution.
  \item mark all elements $T\in \mathcal{N}$
  $$ \mathcal{N}:=\operatorname{argmin}\{ |  \{ \widehat{\mathcal{N}} \subset \mathcal{T}_{h}: \sum_{T\in \widehat{\mathcal{N}}} \Xi^2(T) \geq \theta \sum_{T\in \mathcal{T}_{h}} \Xi^2(T) \}  |  \}$$
  \item estimate local analyticity, i.e.~compute Legendre coefficients of
  $$v_{hp}|_T(\Psi_T(x))=\sum_{j=0}^{p_T} a_{i}L_i(x), \qquad a_{i}=\frac{2i+1}{2} \int_{-1}^1 v_{hp}|_T(\Psi_T(x)) L_i(x) \;dx$$  
  and use a least square approach to compute the slop $m$ of $|\log |a_{i}||=mi+b$, for each direction of $u^{hp}$ on $\Gamma_\Sigma$, of $\psi^{hp}$ on $\Gamma_D$, respectively. If $e^{-m}\leq \delta$ for all directions then $p$-refine, else $h$-refine marked element $E$. If $p_E=0$ always $p$-refine to have a decision basis next time.
  \item refine marked elements based on the decision in 2(e).  
  \end{enumerate}
\end{enumerate}
\hrule
\end{alg}

\section{Numerical experiments}
\label{sec:NumExp}
For the adaptivity algorithm \ref{alg:adaptive} we use D\"orfler marking with marking parameter $\theta=0.3$, and the Legendre expansion strategy \cite{houston2005note} with $\delta=0.5$ for the decision between $h$- and $p$-refinement. The error of the auxiliary problem is estimated by a bubble error estimate as in e.g.~\cite{Banz2015} and the non-localized norms in the delamination related error contributions are approximated by scaled $L^2$-norms, namely $hp^{-1}\| (u^\varepsilon_{hp,n}-g)^+\|^2_{L^2(\Gamma_C)}$ and $h^{-1}p \|( \lambda_{hp}^\varepsilon - DJ_{\varepsilon}(\mathbf{u}^\varepsilon_{hp}))^- \|^2_{L^2(\Gamma_C)}$.
The integral $\langle D J_{\varepsilon}(\mathbf{u}^{\varepsilon}_{hp}), \mathbf{v}_{hp} -\mathbf{u}^{\varepsilon}_{hp} \rangle_{\Gamma_C}$ is computed by a composite Gauss-Quadrature with $\lceil p_T \rceil+17$ quadrature points per element.
If not mentioned otherwise the regularization parameter is $\varepsilon=10^{-4}$. 

\begin{figure}[tbp]
  \centering 
	\includegraphics[trim = 0mm 1mm 14mm 9mm, clip,width=55.0mm, keepaspectratio]{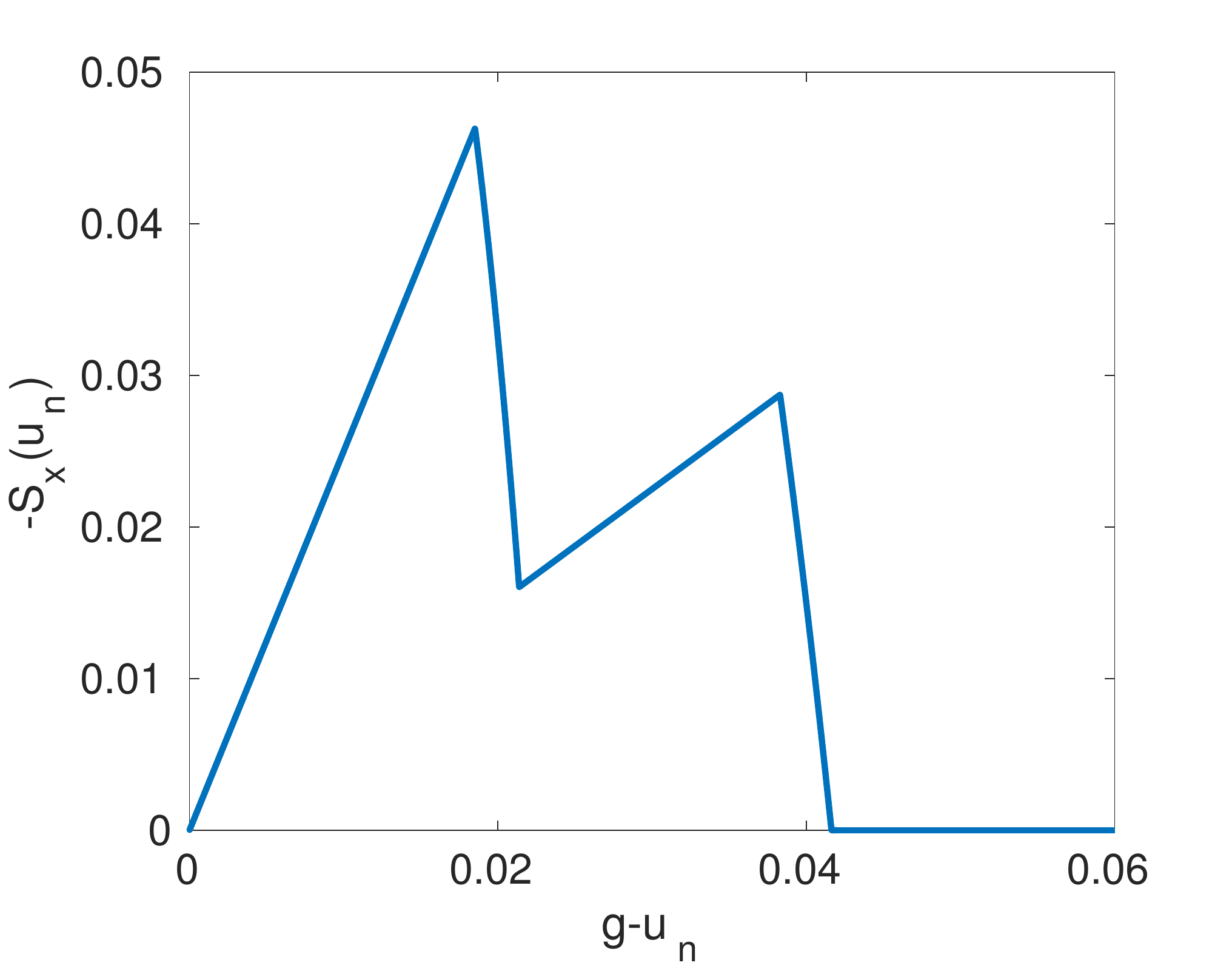}
  \caption{Regularized delamination law $S_x$ for $\varepsilon=10^{-4}$}
  \label{fig:regDelamLaw}
\end{figure}

For the numerical experiments we choose $\Omega=(0,1/2)^2$, $\Gamma_D=\{0\} \times [0,1/2]$, $\Gamma_C=(0,1/2] \times \{0\} $, $\Gamma_N=\partial \Omega \setminus \left( \Gamma_D \cup \Gamma_C \right)$. The material parameters are $E=5$, $\nu=0.45$, $\mathbf{f} \equiv 0$, ${\mathbf t} = 0.25$ on $[1/4,1/2]\times\{1\}$ and zero elsewhere, $g=0$. The delamination law is given via 
\begin{align*}
f(u_n(x))&=\min\{g_1(g(x)-u_n(x)),g_2(g(x)-u_n(x)),g_3(g(x)-u_n(x))\}\\
&= -\max\{-g_1(-u_n(x)),-g_2(-u_n(x)),-g_3(-u_n(x))\} 
\end{align*}
with 
\begin{gather*}
 g_1(y)=\frac{A_1}{2t_1}y^2, \qquad  g_2(y)=b_2(y^2-t_1^2) +d_2, \qquad  g_3(y)=d_3
 \end{gather*}
 and parameters 
 \begin{gather*}
        A_1=0.05, \quad
        A_2=0.03,\quad     
        t_1=0.02, \quad t_2=0.04, \\
        b_2=\frac{A_2}{2t_2},   \quad d_2=A_1\frac{t_1}{2},  \quad d_3=b_2(t_2^2-t_1^2)+d_2 .
 \end{gather*}
The regularized delamination law $S_x$ with regularization parameter $\varepsilon=10^{-4}$ is plotted in Figure~\ref{fig:regDelamLaw}. The characteristic saw tooth shape is already present, but the absolute value in the tips and the slope approximating the jump are still noticeable coarse approximated.

\begin{figure}[tbp]
  \centering \mbox{
  \subfigure[Deformation]{ 
	\includegraphics[trim = 12mm 0mm 12mm 7mm, clip,width=51.0mm, keepaspectratio]{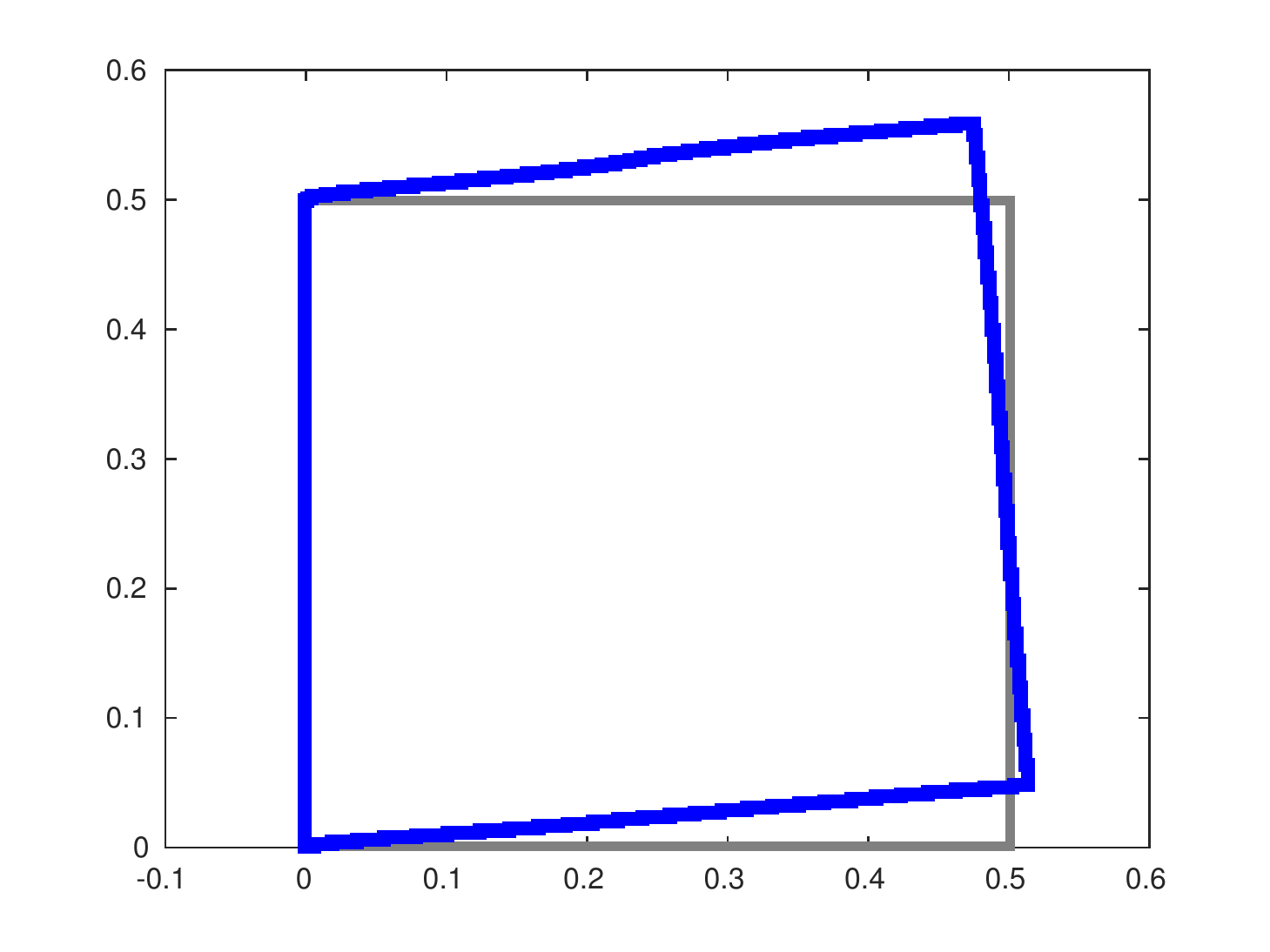}}  \qquad
  \subfigure[Contact stresses]{
	\includegraphics[trim = 2mm 1mm 11mm 6mm, clip,width=55.0mm, keepaspectratio]{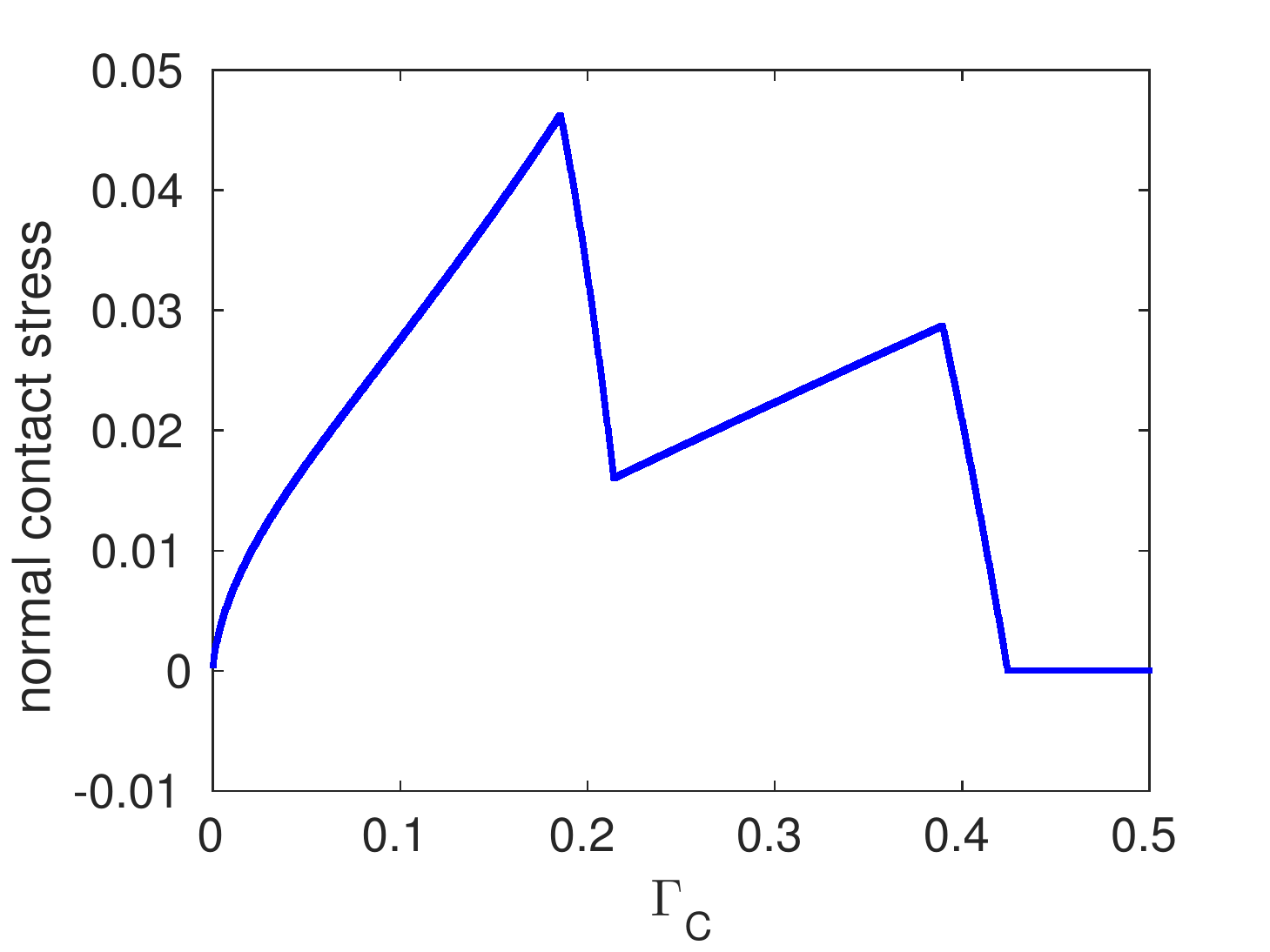}} } 
  \caption{Uniform $h$-version, p=1, 16384 elements with 4096 elements on $\Gamma_C$}
  \label{fig:discreteSol}
\end{figure}

The discrete Lagrange multiplier $\lambda^\varepsilon_{hp}$ is obtain by solving \eqref{eq:reconstruct_lambda_h} where $\psi_i$ are discontinuous, piecewise polynomials on $\Gamma_C$ on a one time coarsened mesh $(H=2h)$ with polynomial degree reduced by one $(q=p-1)$ compared to the mesh and polynomial degree distribution of $\mathbf{u}^\varepsilon_{hp}$. Figure~\ref{fig:discreteSol} displays the deformation of the rectangle and the normal stresses on $\Gamma_C$ obtained from the lowest order uniform $h$-method with 16384 elements and regularization parameter $\varepsilon=10^{-4}$. The normal stress on $\Gamma_C$, Figure~\ref{fig:discreteSol}~(b), reflects the delamination law from Figure~\ref{fig:regDelamLaw} well.

\begin{figure}[tbp]
  \centering
  \includegraphics[trim = 14mm 7mm 15mm 8mm, clip,width=70.0mm, keepaspectratio]{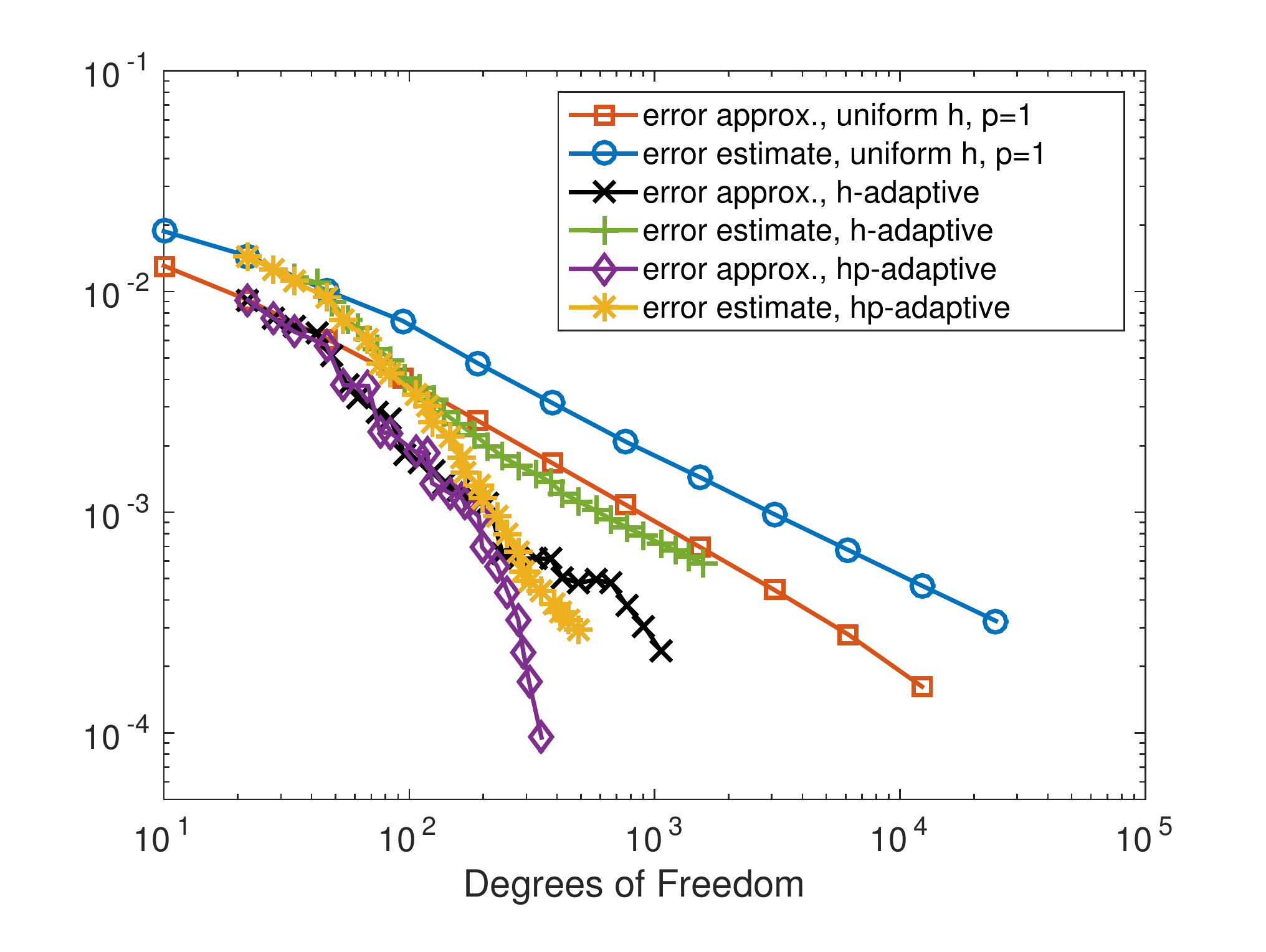}
  \caption{$\sqrt{\|\mathbf{u}_{fine}-\mathbf{u}_{hp}\|_P^2 + \|\lambda_{fine}-\lambda_{hp}\|_V^2}$ and error estimate for different families of discrete solutions, $\varepsilon=10^{-4}$}
  \label{fig:error}
\end{figure}

Figure~\ref{fig:error} displays the reduction of the error in $(\mathbf{u},\lambda)$ and of the error estimate.
Since the exact solution is not known, we compute the error approximately by $\|\mathbf{u}_{fine}-\mathbf{u}_{hp}\|_P$ and $\|\lambda_{fine}-\lambda_{hp}\|_V$, with norms induced by the Poincar\'{e}-Steklov operator, single layer potential acting on $\Gamma_C$, respectively. The pair $(\mathbf{u}_{fine},\lambda_{fine})$ is a very fine (last) approximation for each sequence of discretization. These error expressions are only meaningful if the distance is sufficiently large, which at the end of the different discretization sequences is not true and, therefore, are then omitted. Hence, the error curves are "shorter" then the error estimate curves.\\
The uniform $h$-version with $p=1$ exhibits an experimental order of convergence (eoc) of 0.65 for the error and of 0.53 for the error estimate, all w.r.t.~degrees of freedom (dof). Hence, the error estimate is only reliable but not efficient. For the $h$-adaptive scheme with $p=1$, the eoc is 0.84 for the error and only 0.56 for the error estimate over the last 20, 10 iterations, respectively. In case of $hp$-adaptivity the eoc appears to be exponential for the error and 1.52 for the error estimate also over the last 10 iterations.

\begin{figure}[tbp]
  \centering \mbox{
  \subfigure[uniform $h$-version with $p=1$]{
	\includegraphics[trim = 14mm 7mm 15mm 8mm, clip, clip,width=58.0mm, keepaspectratio]{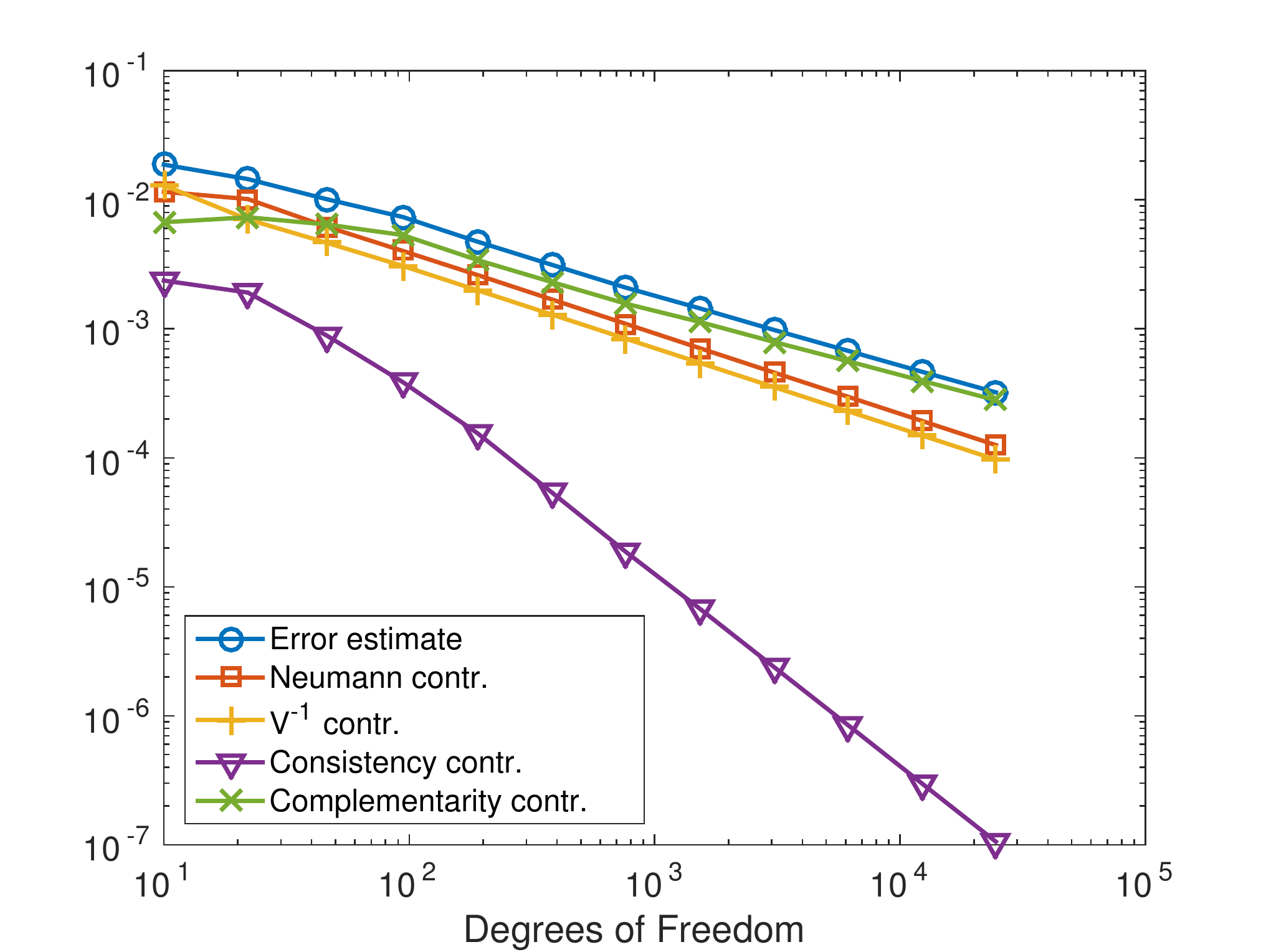}} \ 
  \subfigure[$h$-adaptive with $p=1$ ]{
	\includegraphics[trim = 14mm 7mm 15mm 8mm, clip,width=58.0mm, keepaspectratio]{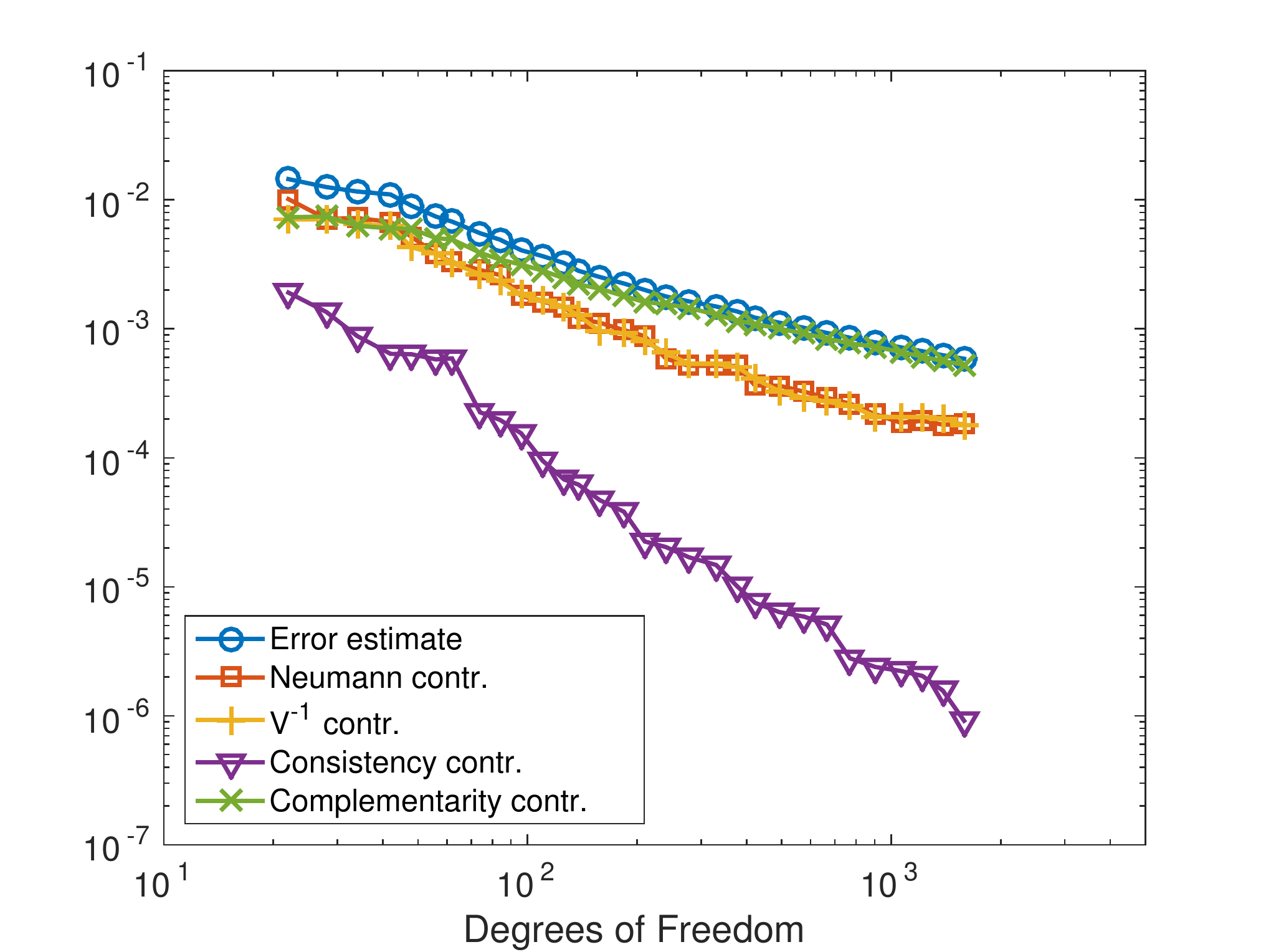}} } \\
	\mbox{
  \subfigure[$hp$-adaptive]{
	\includegraphics[trim = 14mm 7mm 15mm 8mm, clip,width=58.0mm, keepaspectratio]{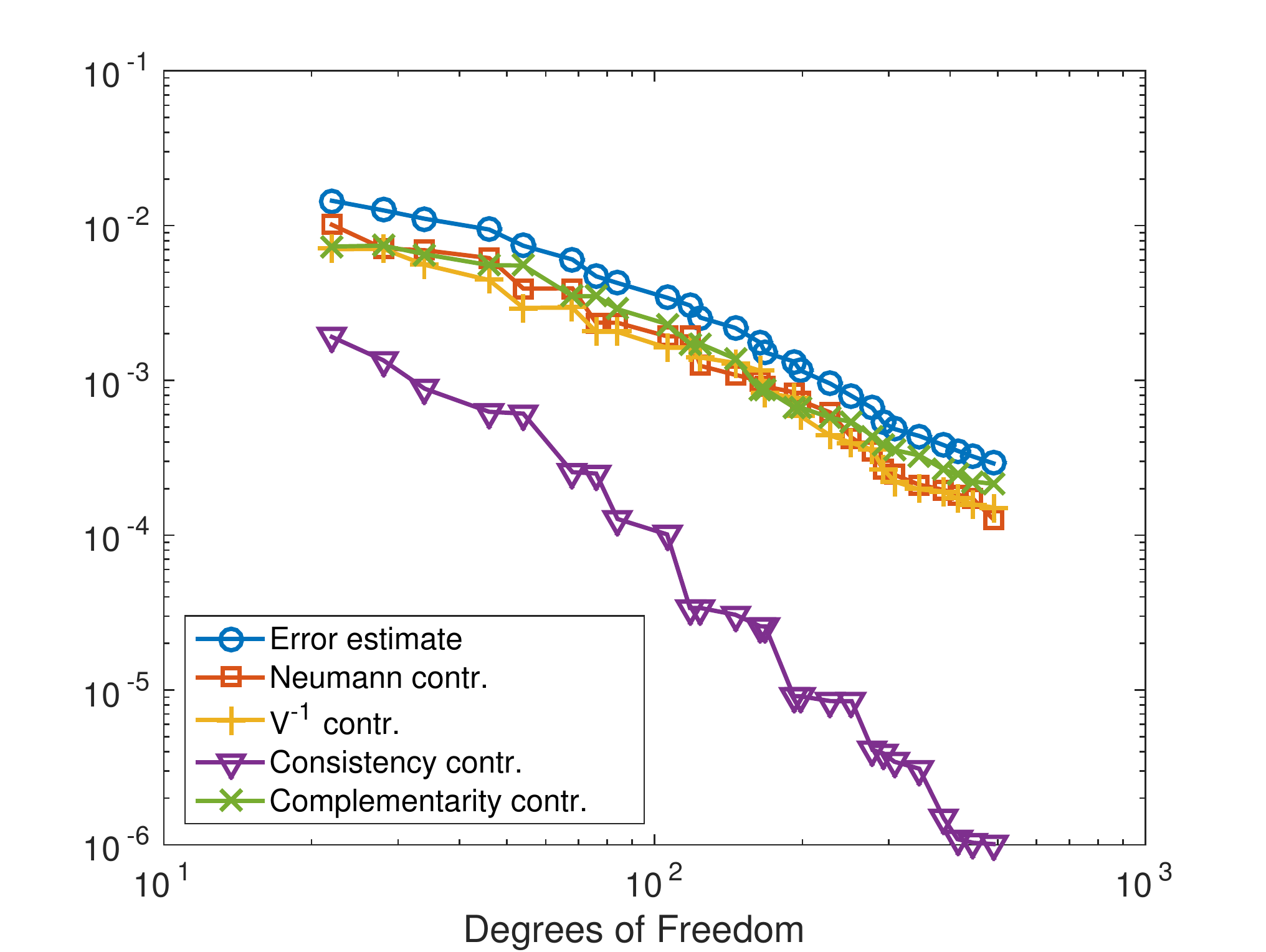}} 
  }
  \caption{Error contributions of the error estimates}
  \label{fig:errorContributions}
\end{figure}

\begin{figure}[tbp]
  \centering \mbox{
  \subfigure[$h$-adaptive, mesh nr.~10 (inner), nr.~20 (outer)]{
	\begin{overpic}	[trim = 28.5mm 10mm 22.5mm 5.5mm, clip, clip,width=55.0mm, keepaspectratio]{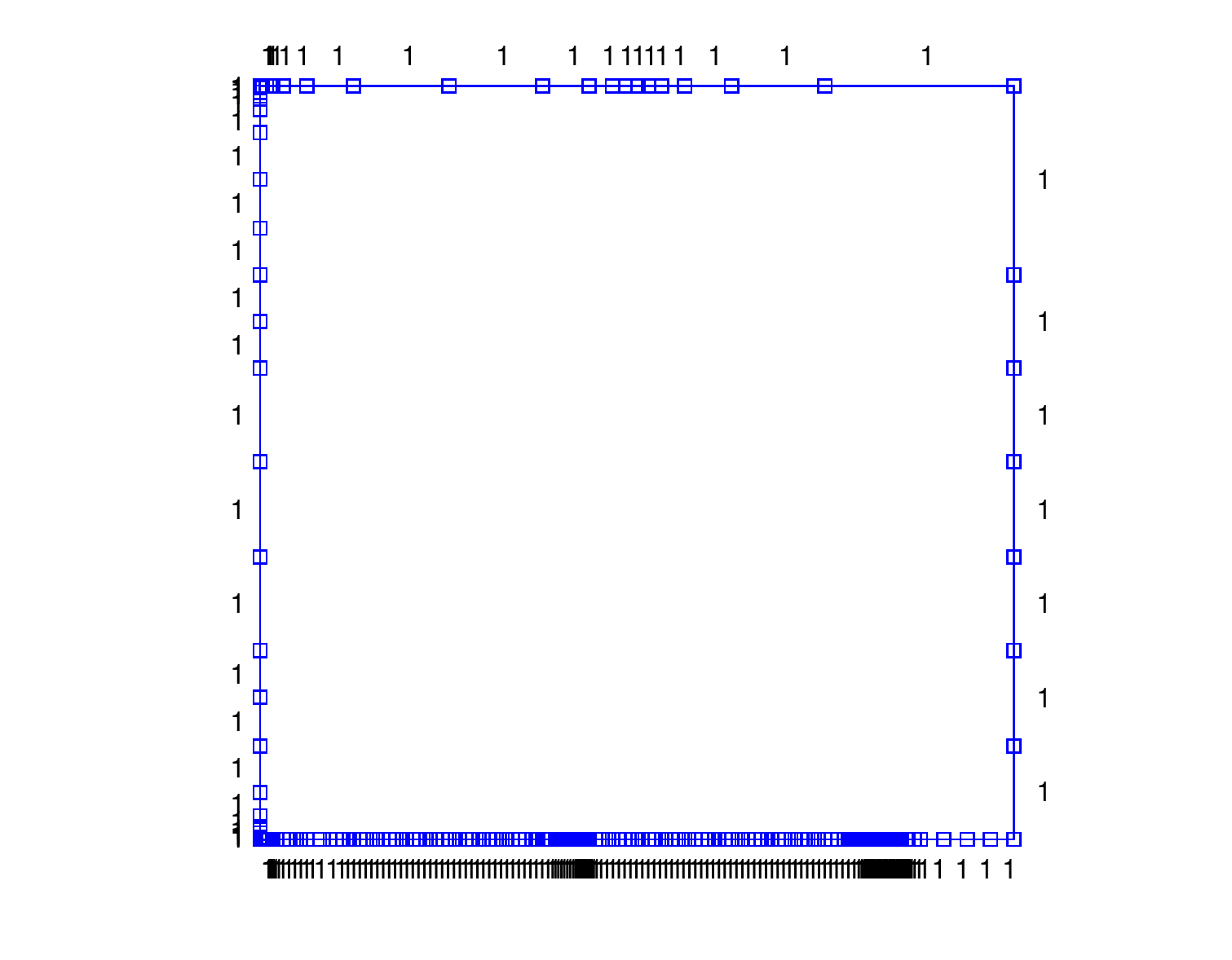}	
	\put(6.5,8){\includegraphics[trim = 29.5mm 7mm 23mm 3mm, clip, clip,width=47.0mm, keepaspectratio]{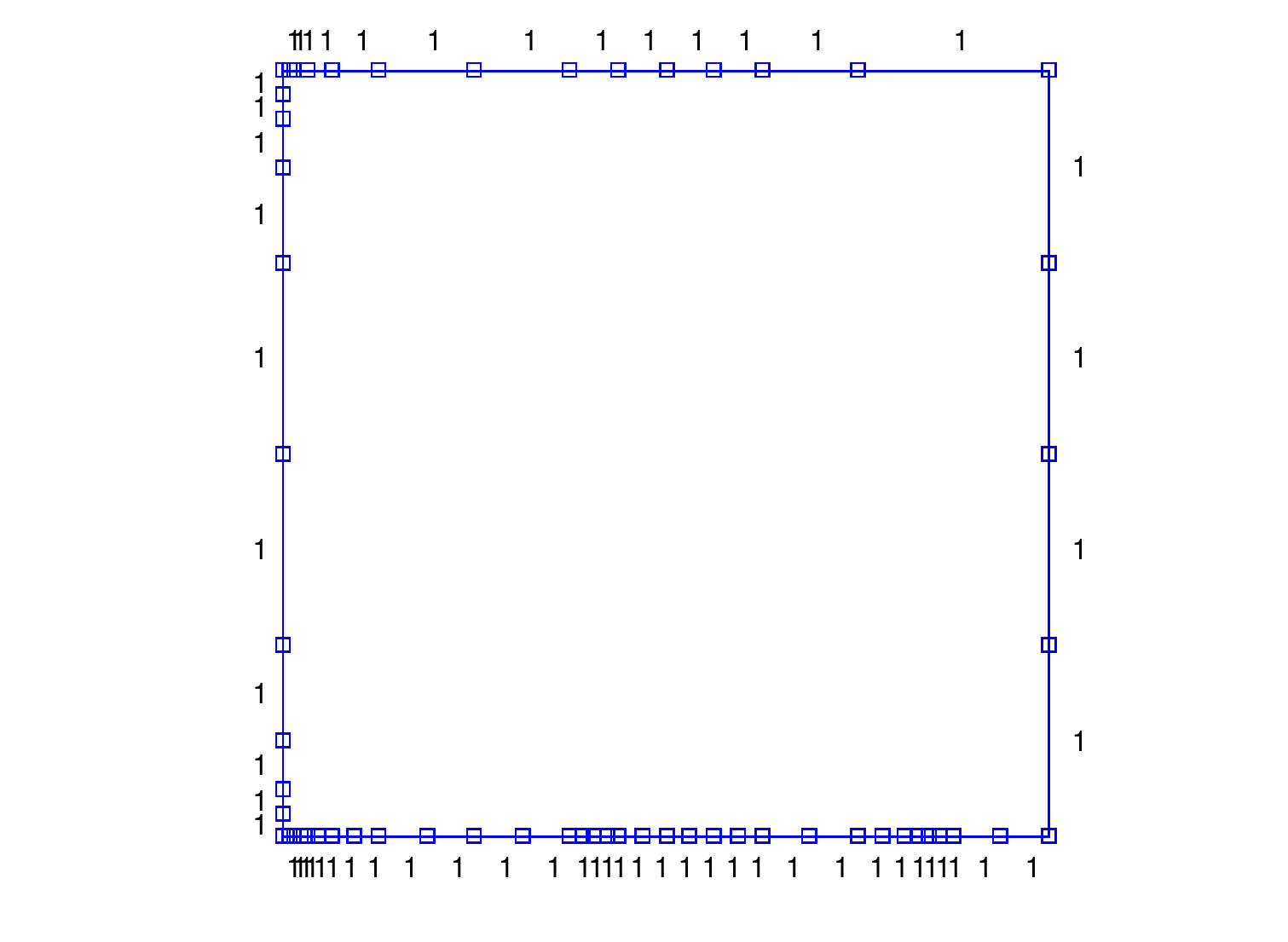}}
	\end{overpic}
}
 \ \ 
  \subfigure[$hp$-adaptive, mesh nr.~15 (inner), nr.~25 (outer)]{
	\begin{overpic}	[trim = 28.5mm 10mm 22.5mm 5.5mm, clip, clip,width=55.0mm, keepaspectratio]{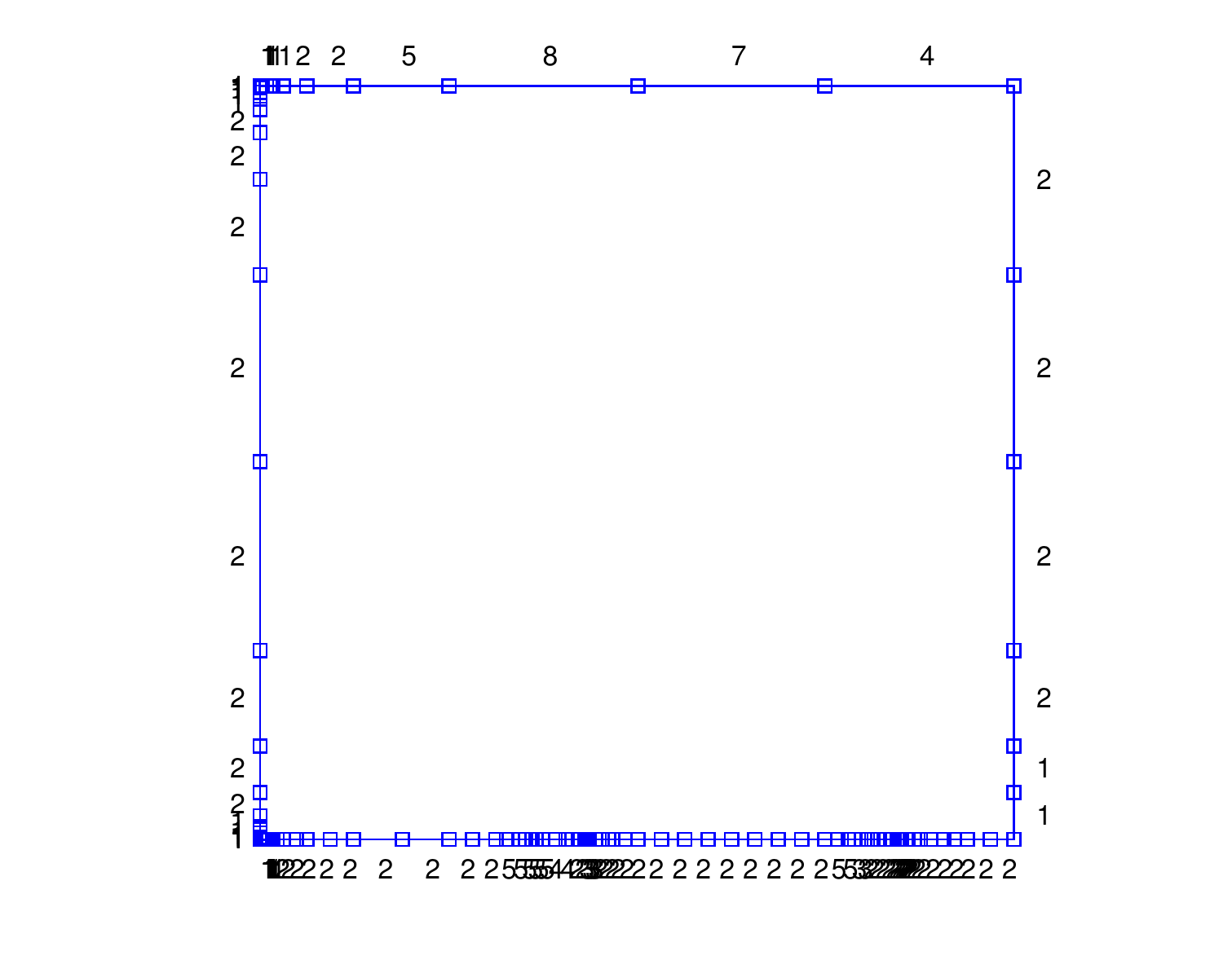}	
	\put(6.7,8){\includegraphics[trim = 27.5mm 7mm 21.5mm 3mm, clip, clip,width=47.0mm, keepaspectratio]{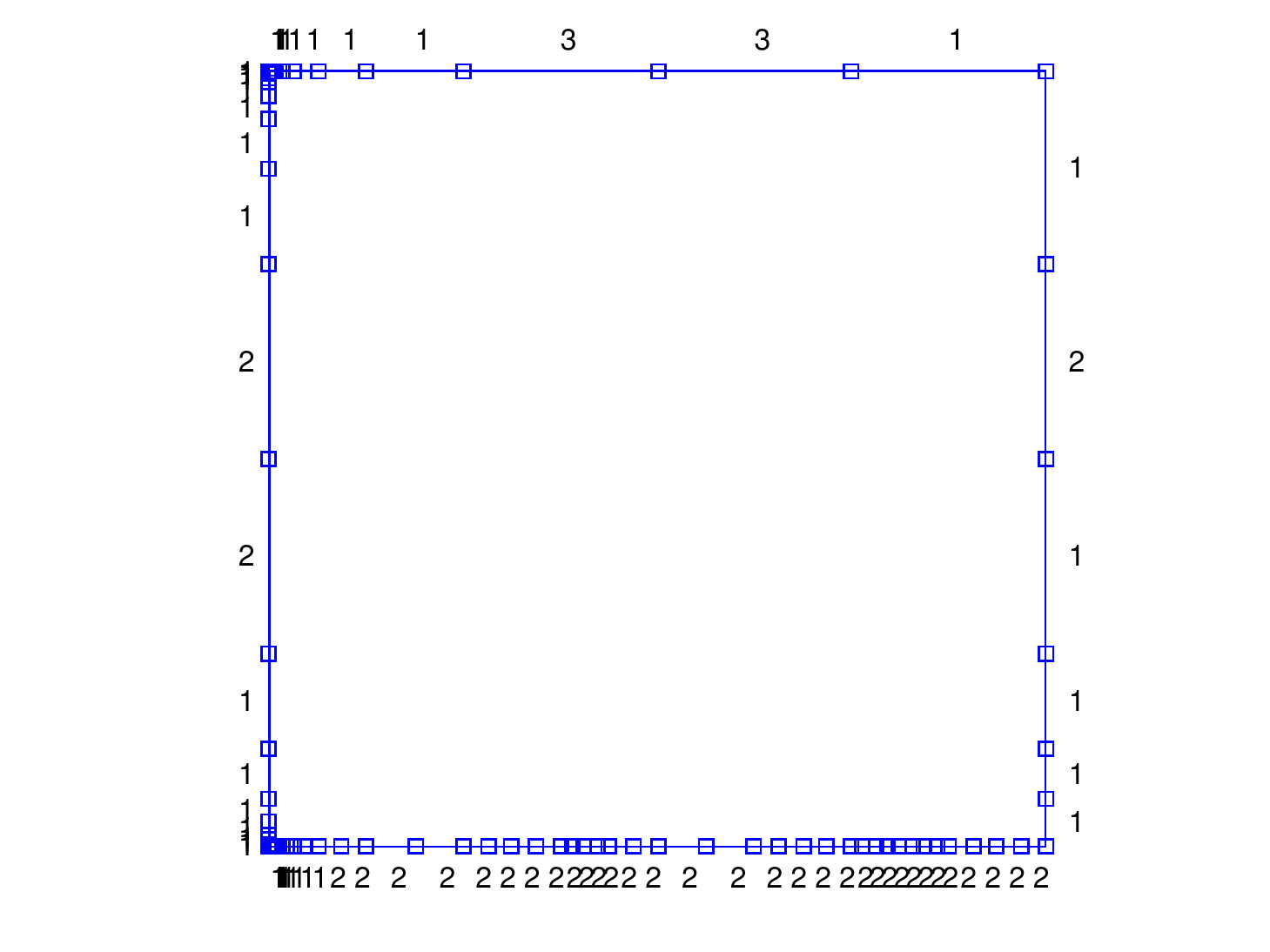}}
	\end{overpic}
} 
}
   \caption{Adaptively generated meshes}
   \label{fig:adaptiveMeshes}
\end{figure}

Looking at the single contributions of the error estimate plotted in Figure~\ref{fig:errorContributions} we find that in all cases the complementarity contribution is the dominant error contribution with the lowest order of convergence and that the consistency contribution has always a significantly higher order of convergence. Given the nature of this numerical experiment, the violation of the non-penetration condition is always zero and, therefore, is omitted in these plots.
For the uniform $h$-version, Figure~\ref{fig:errorContributions}~(a), we have that the Neumann and $V^{-1}$ contributions converge at $0.62$ almost as fast as the error approximation. Given the nature, that $\lambda_{h,1}$ is piecewise constant with double mesh size, $-\mathbf{u}_{h,1} \cdot \mathbf{n}$ is linearly increasing and thus $-S_x(\mathbf{u}_{h,1} \cdot \mathbf{n})$ in the simplest case as well, the complementarity contribution is zero on the left element and non-zero on the right element. This oscillatory case occurs over a large fraction of $\Gamma_C$ and thus the $h$-adaptive scheme performs almost uniform mesh refinements on $\Gamma_C$, c.f.~Figure~\ref{fig:adaptiveMeshes}(a). Hence, the little gain in $h$-adaptivity, Figure~\ref{fig:errorContributions}~(b), for the error estimate. This oscillatory effect can be reduced by increasing the polynomial degree as in $hp$-adaptivity, which balances out the single error contributions, Figure~\ref{fig:errorContributions}~(c) and is thus capable of producing more localized refinements, c.f.~Figure~\ref{fig:adaptiveMeshes}(b). \\

\begin{figure}[tbp]
  \centering 
	\includegraphics[trim = 14mm 7mm 15mm 8mm, clip,width=65.0mm, keepaspectratio]{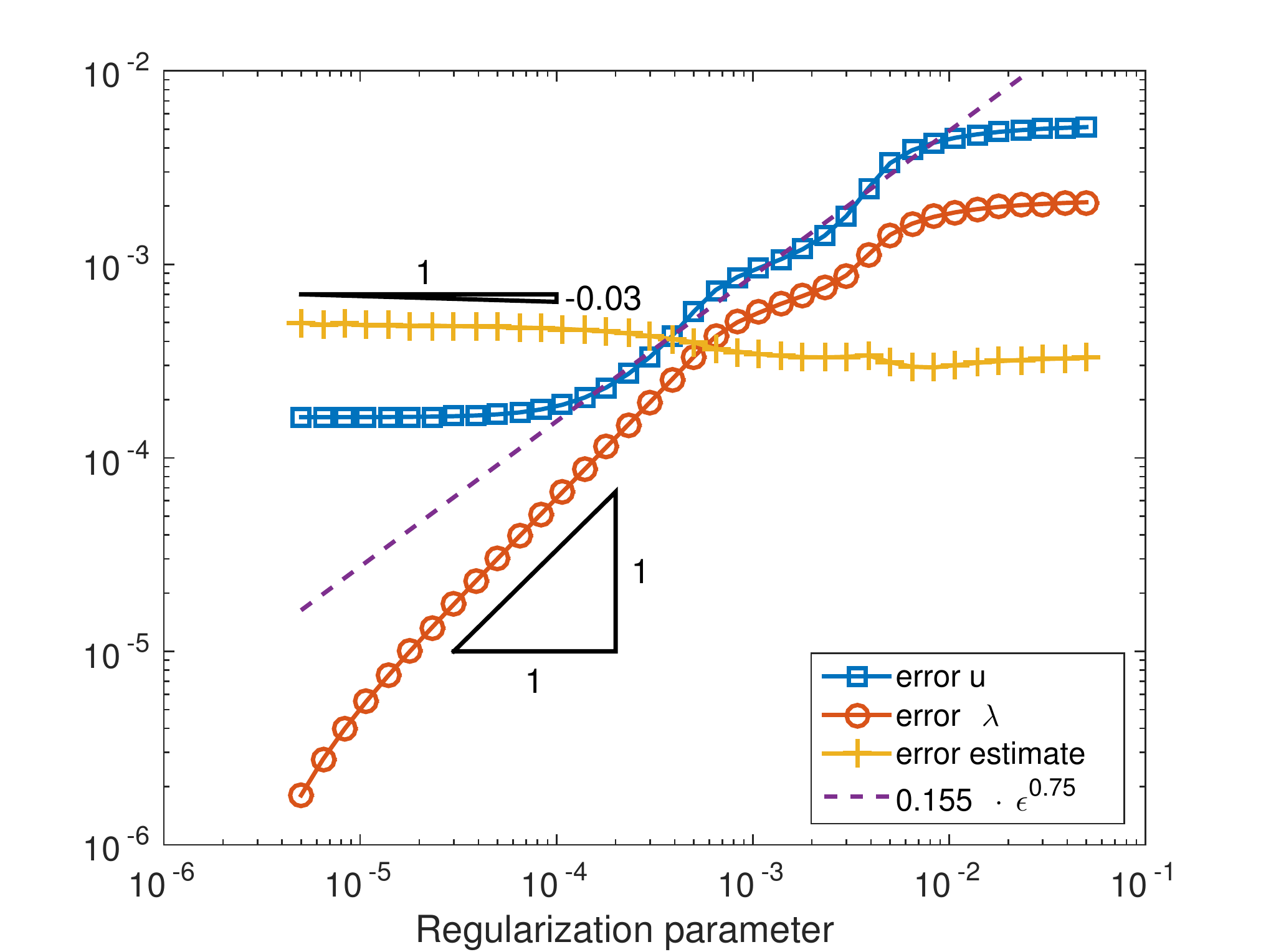}
  \caption{Error and error estimate for different regularization parameters $\varepsilon$, uniform $h$-version 8192 elements, $p=1$}
  \label{fig:errorInfluenceReg}
\end{figure}

Finally, Figure~\ref{fig:errorInfluenceReg} shows the influence of the regularization parameter on the errors $\|\mathbf{u}_{fine}-\mathbf{u}_{h,1}^\varepsilon\|_P$ and $\|\lambda_{fine}-\lambda_{h,1}^\varepsilon\|_V$ where $(\mathbf{u}_{fine},\lambda_{fine})$ is the Galerkin solution to a uniform mesh with 16384 elements, $p=1$ and regularization parameter $\varepsilon=2.5 \cdot 10^{-6}$. Here, the pair $(\mathbf{u}_{h,1}^\varepsilon,\lambda_{h,1}^\varepsilon)$ is computed on a uniform mesh with 8192 elements and the regularization parameter is varied. We find that for very large $\varepsilon$, the regularized delamination law has lost its saw tooth characteristic and small variations of $\varepsilon$ have no noticeable influence. For very small $\varepsilon$ the discretization error is dominant which sets in for $\mathbf{u}$ much earlier than for $\lambda$. In the intermediate range the eoc w.r.t.~$\varepsilon$ is 0.75 for $\mathbf{u}$, whereas $\lambda$ displays an (asymptotic) eoc of 1.0 w.r.t.~$\varepsilon$. The error estimate for the regularized problem increases in $\varepsilon$ slightly.\\

{\bf Acknowledgment} The authors are grateful to J. Gwinner  for the helpful discussions and comments.

\bibliographystyle{spmpsci}      
\bibliography{bib_HVI}   


\end{document}